\documentclass[aop,MSNbibl,seceqn,citesort,dvips]{arximspdf}

\doi{10.1214/11-AOP649}
\volume{40}
\issue{3}
\pubyear{2012}
\firstpage{1041}
\lastpage{1068}

\makeatletter
\newcommand{\fraca}[2]{{#1}/{#2}}

\newtheorem{theorem}{Theorem}[section]
\newproclaim{definition}[theorem]{Definition}
\newtheorem{lemma}[theorem]{Lemma}
\newtheorem{proposition}[theorem]{Proposition}
\newproclaim{remark}[theorem]{Remark}
\makeatother

\begin{document}
\begin{frontmatter}

\title{Feynman--Kac formula for the heat equation driven by fractional noise
with Hurst parameter $\bolds{H<1/2}$}
\runtitle{Feynman--Kac formula}

\begin{aug}
\author[A]{\fnms{Yaozhong} \snm{Hu}\corref{}\ead[label=e1]{hu@math.ku.edu}},
\author[A]{\fnms{Fei} \snm{Lu}\ead[label=e2]{feilu@math.ku.edu}}
\and
\author[A]{\fnms{David} \snm{Nualart}\thanksref{t2}\ead[label=e3]{nualart@math.ku.edu}}
\runauthor{Y. Hu, F. Lu and D. Nualart}
\affiliation{University of Kansas}
\address[A]{Department of Mathematics \\
University of Kansas \\
Lawrence, Kansas 66045 \\
USA\\
\printead{e1}\\
\phantom{\textsc{E-mail:}\ }\printead*{e2}\\
\phantom{\textsc{E-mail:}\ }\printead*{e3}} 
\end{aug}
\thankstext{t2}{Supported in part by the NSF Grant DMS-09-04538.}

\received{\smonth{7} \syear{2010}}
\revised{\smonth{1} \syear{2011}}

%
\begin{abstract}
In this paper, a Feynman--Kac formula is established for stochastic partial
differential equation driven by Gaussian noise which is, with respect to
time, a fractional Brownian motion with Hurst parameter $H<1/2$. To
establish such a formula, we introduce and study a nonlinear stochastic
integral from the given Gaussian noise. To show the Feynman--Kac integral
exists, one still needs to show the exponential integrability of nonlinear
stochastic integral. Then, the approach of approximation with techniques
from Malliavin calculus is used to show that the Feynman--Kac integral
is the
weak solution to the stochastic partial differential equation.
\end{abstract}

%

\end{frontmatter}

\section{Introduction}\label{sec1}

Consider the stochastic heat equation on $\mathbb{R}^{d}$%
%
\begin{equation}
\cases{\displaystyle
\frac{\partial u}{\partial t}=\frac{1}{2}\Delta u+u\,\frac{\partial W}{
\partial t}(t,x) ,&\quad $t\geq0 , x\in\mathbb{R}^{d} $,\vspace*{1pt}\cr\displaystyle
u(0,x)=u_{0}(x) ,
}
\label{she}
\end{equation}
where $u_{0}$ is a bounded measurable function and $W=\{W(t,x), t\geq
0,x\in
 \mathbb{R}^{d}\}$ is a fractional Brownian motion of Hurst
parameter $H
\in ( \frac{1}{4},\frac{1}{2} ) $ in time and it has a spatial
covariance $Q(x,y)$, which is locally $\gamma$-H\"{o}lder continuous (see
Section~\ref{sec2} for precise meaning of this condition), with $\gamma>2-4H$. We
shall show that the solution to (\ref{she}) is given by
%
\begin{equation}
u(t,x)=E^{B} \biggl[ u_{0}(B_{t}^{x})\exp\int
_{0}^{t}W(ds,B_{t-s}^{x}) \biggr]
 , \label{e.1.3}
\end{equation}
where $B=\{B_{t}^{x}=B_{t}+x,t\geq0,x\in\mathbb{R}^{d}\}$ is a $d$-dimensional Brownian motion starting at $x\in\mathbb{R}^{d}$, independent
of $W$.

This is a generalization of the well-known Feyman--Kac formula to the
case of
a random potential of the form  $\frac{\partial W}{\partial t}(t,x)$.
Notice that the integral $\int_{0}^{t}W(ds,B_{t-s}^{x})$ is a nonlinear
stochastic integral with respect to the fractional noise $W$. This type of
Feynman--Kac formula was mentioned as a~conjecture by Mocioalca and
Viens in
\cite{M-V05}.

There exists an extensive literature devoted to Feynman--Kac formulas for
stochastic partial differential equations. Different versions of the
Feynman--Kac formula have been established for a variety of random
potentials. See, for instance, a Feynman--Kac formula for anticipating SPDE
proved by Ocone and Pardoux \cite{O-P93}. Ouerdiane and Silva \cite{OS02}
give a generalized Feynman--Kac formula with a convolution potential by
introducing a generalized function space. Feynman--Kac formulas for L\'{e}vy
processes are presented by Nualart and Schoutens \cite{Nu-Sch01}.

However, only recently a Feynman--Kac formula has been established by Hu
 et al.  \cite{hunuso} for random potentials associated with the
fractional Brownian motion. The authors consider the following stochastic
heat equation driven by fractional noise
%
\begin{equation}
\cases{\displaystyle
\frac{\partial u}{\partial t}=\frac{1}{2}\Delta u+u\,\frac{\partial
^{d+1}W}{%
\partial t\,\partial x_{1}\cdots\partial x_{d}}(t,x) ,&\quad $t\geq0,
x\in
\mathbb{R}^{d} $,\cr\displaystyle
u(0,x)=u_{0}(x) ,%
}
\label{e.1.1}
\end{equation}
where $W=\{W(t,x), t\geq0,x\in \mathbb{R}^{d}\}$ is fractional Brownian
sheet with Hurst parameter $(H_{0},H_{1},\ldots ,H_{d})$. They show
(\cite{hunuso}, Theorem 4.3) that if $H_{1},\ldots,\allowbreak H_{d}\in(\frac{1}{2},1) $,
and $2H_{0}+H_{1}+\cdots+H_{d}>d+1$, then the solution $u(t,x)$ to the
above stochastic heat equation is given by
%
\begin{equation}
u(t,x)=E^{B} \biggl[ f(B_{t}^{x})\exp \biggl( \int_{0}^{t}\int
_{\mathbb{R}%
^{d}}\delta(B_{t-r}^{x}-y)W(dr,dy) \biggr)  \biggr]  , \label{feynman}
\end{equation}
where $B=\{B_{t}^{x}=B_{t}+x,t\geq0,x\in\mathbb{R}^{d}\}$ is a $d$-dimensional Brownian motion starting at $x\in\mathbb{R}^{d}$, independent
of $W$. The condition $2H_{0}+H_{1}+\cdots+H_{d}>d+1$ is shown to be sharp
in that framework. Since the $H_{i}$, $i=1,\ldots ,d$, cannot take value
greater or equal to $1$, this condition implies that $H_{0}>\frac{1}{2}$.

We remark that if $B^{H_{0}}=\{B_{t}^{H_{0}},t\geq0\}$ is a fractional
Brownian motion with Hurst parameter $H_{0}>\frac{1}{2}$, then the
stochastic integral $\int_{0}^{T}f(t)\,dB_{t}^{H_{0}}$ is well defined
for a
suitable class of distributions $f$, and in this sense the above
integral $%
\int_{0}^{t}\int_{\mathbb{R}^{d}}\delta(B_{t-r}^{x}-y)W(dr,dy)$ is well
defined for any trajectory of the Brownian motion $B$.  If
$H_{0}<\frac{1%
}{2}$,  this is no longer true and we can integrate only functions
satisfying some regularity conditions. For this reason, it is not possible
to write a Feynman--Kac formula for the equation (\ref{e.1.1}) with
$H_{0}<%
\frac{1}{2}$.

Notice that for $d=1$ and $H_{0}=H_{1}=\frac{1}{2}$ (space--time white
noise) a Feyn\-man--Kac formula can not be written for equation (\ref{e.1.1}),
but this equation has a~unique mild solution when the stochastic
integral is
interpreted in the It\^{o} sense. A renormalized Feynman--Kac formula with
Wick exponential has been obtained in this case by Bertinin and Cancrini
\cite{BC95}. More generally, if the product appearing in (\ref
{e.1.1}) is
replaced by Wick product, Hu and Nualart~\cite{H-N07} showed that a formal
solution can be obtained using chaos expansions.

In the present paper, we are concerned with the case $H_{0}<\frac
{1}{2}$, but
we use a random potential of the form $\frac{\partial W}{\partial t}(t,x)$.
One of the main obstacles to overcome is to define the stochastic
integral $%
\int_{0}^{t}W(ds,B_{t-s}^{x})$. We start with the construction of a general
nonlinear stochastic integral $\int_{0}^{t}W(ds,\phi_{s})$ where
$\phi$ is
a H\"{o}lder continuous function of order $\alpha>\frac{1}{\gamma}(1-2H)$.
It turns out that the irregularity in time of $W(t,x)$ is compensated
by the
above H\"{o}lder continuity of $\phi$ through the covariance in space, with
an appropriate application of the fractional integration by parts technique.
Let us point out that $\int_{0}^{t}W(ds,\phi_{s})$ is well defined
for all H\"{o}lder continuous function $\phi$ with $\alpha>\frac{1}{\gamma
}(\frac{1%
}{2}-H)$, and we consider here only the case $\alpha>\frac{1}{\gamma}
(1-2H) $ because this condition is required when we show that $u (
t,x ) $ is a weak solution to~(\ref{she}). Furthermore, the
condition $%
\alpha>\frac{1}{\gamma}(1-2H)$ also allows us to obtain an explicit
formula for the variance of $\int_{0}^{t}W(ds,\phi_{s})$. Contrary to
\cite{hunuso}, it is rather simpler to show that $\int_{0}^{t}W(ds,B_{t-s}^{x})$
is exponentially integrable. A~by-product is that $u(t,x)$ defined by %
(\ref{e.1.3}) is almost surely H\"{o}lder continuous of order which
can be arbitrarily close to $H-\frac{1}{2}+\frac{\gamma}{4} $ from below.
Let us also mention recent work on stochastic integral \cite{HJT} and
\cite{KR} with general Gaussian processes which can be applied to the
case $H<\frac{1}{2}$.

Another main effort of this paper is to show that $u ( t,x ) $
defined by (\ref{e.1.3}) is a solution to (\ref{she}) in a weak sense (see
Definition \ref{def weak solu}). As in \cite{hunuso}, this is done by using
an approximation scheme together with techniques of Malliavin calculus. Let
us point out that in the definition of $\int_{0}^{t}W(ds,\phi_{s})$
one can
use a~one-side approximation, but it is necessary to use symmetric
approximations (as well as the condition $H>\frac{1}{2}-\frac{\gamma
}{4} $%
) to show the convergence of the trace term (\ref{Trace}).

We also discuss the corresponding Skorohod-type equation, which corresponds
to taking the Wick product in \cite{H-N07}. We show that a unique mild
solution exists for $H\in(\frac{1}{2}-\frac{\gamma}{4} ,\frac
{1}{2})$.%

The paper is organized as follows. Section~\ref{sec2} contains some
preliminaries on
the fractional noise $W$ and some results on fractional calculus which is
needed in the paper. We also list all the assumptions that we make for the
noise $W$ in this section. In Section~\ref{sec3}, we study the nonlinear stochastic
integral appeared in equation (\ref{e.1.3}) by using smooth approximation
and we derive some basic properties of this integral. Section~\ref{sec4}
verifies the
integrability and H\"{o}lder continuity of $u ( t,x ) $.
Section~\ref{sec5}
is devoted to show that $u(t,x)$ is a solution to (\ref{she}) in a weak
sense. Section~\ref{sec6} gives a solution to the Skorohod type equation. The last
section is the \hyperref[appm]{Appendix} with some technical results used along the paper.

\section{Preliminaries}\label{sec2}

Fix $H\in(0,\frac{1}{2})$ and denote by $R_{H}(t,s)=\frac{1}{2} (
t^{2H}+s^{2H}-|t-s|^{2H} ) $ the covariance function of the fractional
Brownian motion of Hurst parameter $H$. Suppose that $W=\{W(t,x),t\geq
0,x\in\mathbb{R}^{d}\}$ is a mean zero Gaussian random field, defined
on a
probability space $({\Omega},\mathcal{F},P)$, whose covariance
function is
given by
\[
E(W(t,x)W(s,y))=R_{H}(t,s)Q(x,y),
\]
where $Q(x,y)$ satisfies the following properties for some $M<2$ and
$\gamma
\in(0,1]$:

\begin{longlist}[(Q2)]
\item[(Q1)]  $Q$ is locally bounded: there exists a constant $C_{0}>0$
such that for any $K>0$
\[
Q(x,y)\leq C_{0} ( 1+K ) ^{M}
\]
for any $x,y\in\mathbb{R}^{d}$ such that $ \vert x \vert
, \vert y \vert\leq K$.

\item[(Q2)]   $Q$ is locally $\gamma$-H\"{o}lder continuous: there
exists a constant $C_{1}>0 $ such that for any $K>0$
\[
 \vert Q(x,y)-Q(u,v) \vert\leq C_{1} ( 1+K )
^{M} (
 \vert x-u \vert^{\gamma}+ \vert y-v \vert
^{\gamma
} )
\]
for any $x,y,u,v\in\mathbb{R}^{d}$ such that  $ \vert x
\vert
, \vert y \vert, |u|, |v| \leq K$.
\end{longlist}

Denote by $\mathcal{E}$ the vector space of all step functions on $[0,T]$.
On this vector space~$\mathcal{E}$, we introduce the following scalar product
\[
\bigl\langle\mathbf{1}_{[0,t]},\mathbf{1}_{[0,s]}\bigr\rangle_{\mathcal{H}%
_{0}}=R_{H}(t,s).
\]
Let $\mathcal{H}_{0}$ be the closure of $\mathcal{E}$ with respect to the
above scalar product. Denote by $C^{\alpha}( [ a,b ] )$ the
set of all functions which
is H\"{o}lder continuous of order ${\alpha}$, and denote by $
\|\cdot \|_{\alpha}$ the $\alpha
$-H\"{o}%
lder norm. It is well known that $C^{\alpha} ( [0,T] )
\subset\mathcal{H}_{0}$ for $\alpha>\frac{1}{2}-H$.

Let $\mathcal{H}$ be the Hilbert space defined by the completion of the
linear span of indicator functions $\mathbf{1}_{[0,t]\times\lbrack
0,x]}$, $%
t\in\lbrack0,T]$, $x\in\mathbb{R}^{d}$ under the scalar product%
\[
\bigl\langle\mathbf{1}_{[0,t]\times\lbrack0,x]},\mathbf{1}_{[0,s]\times
\lbrack0,y]}\bigr\rangle_{\mathcal{H}}=R_{H}(t,s)Q(x,y).
\]
In the above formula, if $x_{i}<0$ we assume by convention that
$\mathbf{1}%
_{[0,x_{i}]}=-\mathbf{1}_{[-x_{i},0]}$. The mapping $W\dvtx \mathbf{1}%
_{[0,t]\times\lbrack0,x]}\rightarrow W(t,x)$ can be extended to a linear
isometry between $\mathcal{H}$ and the Gaussian space spanned by $W$.
Then, $%
\{W(h), h\in\mathcal{H}\}$ is an isonormal Gaussian process.

Let $\mathcal{S}$ be the space of random variables $F$ of the form:
\[
F=f(W(\varphi_{1}),\ldots ,W(\varphi_{n})),
\]
where $\varphi_{i}\in\mathcal{H}$, $f\in C^{\infty}(\mathbb
{R}^{n})$, $f$
and all its partial derivatives have polynomial growth. The Malliavin
derivative\vadjust{\goodbreak} $DF$ of an element $F$ in $\mathcal{S}$ is defined as an $%
\mathcal{H}$-valued random variable given by
\[
DF=\sum_{i=1}^{n}\frac{\partial f}{\partial x_{i}}(W(\varphi
_{1}),\ldots ,W(\varphi_{n}))\varphi_{i}.
\]
The operator $D$ is closable from $L^{2}(\Omega)$ into $L^{2}(\Omega,%
\mathcal{H})$ and we define the Sobolev space $\mathbb{D}^{1,2}$ as the
closure of $\mathcal{S}$ with respect to the following norm:
\[
 \| DF \|_{1,2}=\sqrt{%
E(F^{2})+E( \| DF \|
_{\mathcal{H}%
}^{2})}.
\]
The divergence operator $\delta$ is the adjoint of the derivative
operator $%
D$, determined by the duality relationship
\[
E(\delta(u)F)=E(\langle DF,u\rangle_{\mathcal{H}}) \qquad
\mbox{for any
$F\in\mathbb{D}^{1,2}$}.
\]
$\delta(u)$ is also called the Skorohod integral of $u$. We refer to
Nualart \cite{Nu06} for a~detailed account on the Malliavin calculus. For
any random variable \mbox{$F\in\mathbb{D}^{1,2}$} and $\phi\in\mathcal{H},$
%
\begin{equation}
FW ( \phi ) =\delta(F\phi)+ \langle DF,\phi
\rangle_{%
\mathcal{H}}. \label{div_fmla}
\end{equation}

Since we deal with the case of Hurst parameter $H\in(0, 1/2)$, we
shall use
intensively the fractional calculus. We recall some basic definitions and
properties. For a detailed account, we refer to \cite{SKM}.

Let $a,b\in\mathbb{R}$, $a<b$. Let $f\in L^{1} ( a,b ) $
and $%
\alpha>0.$ The left and right-sided fractional integral of $f$ of
order $%
\alpha$ are defined for $x\in ( a,b ) $, respectively, as
\[
I_{a+}^{\alpha}f ( x ) =\frac{1}{\Gamma ( \alpha
 ) }%
\int_{a}^{x} ( x-y ) ^{\alpha-1}f ( y )\,dy
\]
and
\[
I_{b-}^{\alpha}f ( x ) =\frac{1}{\Gamma ( \alpha
 ) }%
\int_{x}^{b} ( y-x ) ^{\alpha-1}f ( y )\,dy .
\]
Let $I_{a+}^{\alpha} ( L^{p} )$ [resp., $I_{b-}^{\alpha
} (
L^{p} ) $] the image of $L^{p} ( a,b ) $ by the
operator $%
I_{a+}^{\alpha}$\break (resp.,~$I_{b-}^{\alpha}$).

If $f\in I_{a+}^{\alpha} ( L^{p} ) $ [resp., $I_{b-}^{\alpha
} (
L^{p} ) $] and $0<\alpha<1$ then the left and right-sided fractional
derivatives are defined by
\begin{equation}
D_{a+}^{\alpha}f ( x ) =\frac{1}{\Gamma ( 1-\alpha
 ) }%
 \biggl( \frac{f ( x ) }{ ( x-a ) ^{\alpha
}}+\alpha
\int_{a}^{x}\frac{f ( x ) -f ( y ) }{ (
x-y )
^{\alpha+1}}\,dy \biggr)   \label{frDl}
\end{equation}
and%
%
\begin{equation}
D_{b-}^{\alpha}f ( x ) =\frac{(-1)^{\alpha}}{\Gamma (
1-\alpha ) } \biggl( \frac{f ( x ) }{ (
b-x ) ^{\alpha
}}+\alpha\int_{x}^{b}\frac{f ( x ) -f ( y )
}{ (
y-x ) ^{\alpha+1}}\,dy \biggr) \label{frDr}
\end{equation}
for all $x\in ( a,b ) $ [the convergence of the integrals
at the
singularity $y=x$ holds point-wise for almost all $x\in (
a,b ) $
if $p=1$ and moreover in $L^{p}$-sense if $1<p<\infty$].

It is easy to check that if $f\in I_{a+ ( b- ) }^{1} (
L^{1} ) ,$
%
\begin{equation}
D_{a+}^{\alpha}D_{a+}^{1-\alpha}f=Df ,\qquad D_{b-}^{\alpha
}D_{b-}^{1-\alpha}f=Df \label{fracint1}
\end{equation}
and
%
\begin{equation}
 ( -1 ) ^{\alpha}\int_{a}^{b}D_{a+}^{\alpha}f (
x )
g ( x )\,dx=\int_{a}^{b}f ( x ) D_{b-}^{\alpha
}g (
x )\,dx \label{fracint2}
\end{equation}
provided that $0\leq\alpha\leq1$, $f\in I_{a+}^{\alpha} (
L^{p} ) $ and $g\in I_{b-}^{\alpha} ( L^{q} )
$ with $%
p\geq1,q\geq1,\frac{1}{p}+\frac{1}{q}\leq1+\alpha.$

It is clear that $D^{\alpha}f$ exists for all $f\in C^{\beta
}([a,b])$ if ${\alpha}<{\beta}$. The following proposition was proved in
\cite{Zahle98}.

\begin{proposition}
\label{p1} Suppose that $f\in C^{\lambda}( [ a,b ] )$ and
$g\in
C^{\mu}( [ a,b ] )$ with $\lambda+\mu>1$. Let ${\lambda
}>\alpha
$ and $\mu>1-\alpha$. Then the Riemann--Stieltjes integral $\int_{a}^{b}fdg$
exists and it can be expressed as%
%
\begin{equation}
\int_{a}^{b}fdg=(-1)^{\alpha}\int_{a}^{b}D_{a+}^{\alpha}f (
t )
D_{b-}^{1-\alpha}g_{b-} ( t )\,dt, \label{e.2.6}
\end{equation}
where $g_{b-} ( t ) =g ( t ) -g ( b ) $.
\end{proposition}

\section{Nonlinear stochastic integral}\label{sec3}

In this section, we introduce the nonlinear stochastic integral that appears
in the Feynman--Kac formula  (\ref{e.1.3})  and obtain some properties
of this integral which are useful in the following sections. The main idea
to define this integral is to use an appropriate approximation scheme. In
order to introduce our approximation, we need to extend the fractional
Brownian field to $t<0$. This can be done by defining $W=\{W(t,x),t\in
\mathbb{R},x\in\mathbb{R}^{d}\}$ as a mean zero Gaussian process
with the
following covariance
\[
E [ W ( t,x ) W ( s,y )  ] =\tfrac
{1}{2} (
 | t | ^{2H}+ | s | ^{2H}-|t-s|^{2H} )
Q (
x,y )  .
\]

For any $\varepsilon>0$, we introduce the following approximation of $%
W(t,x) $:
%
\begin{equation}
W^{\varepsilon}(t,x)= \int_{0}^{t}\dot{W}^{\varepsilon}(s,x)\,ds ,
\label{W_approx}
\end{equation}
where $\dot{W}^{\varepsilon}(s,x)=\frac{1}{2\varepsilon} (
W(s+\varepsilon,x)-W(s-\varepsilon,x) ) $.

\begin{definition}
Given a continuous function $\phi$ on $[0,T]$, define%
\[
\int_{0}^{t}W(ds,\phi_{s})=\lim_{\varepsilon\rightarrow0}\int
_{0}^{t}\dot{%
W}^{\varepsilon}(s,\phi_{s})\,ds,
\]
if the limit exits in $L^{2}(\Omega)$.
\end{definition}

Now we want to find conditions on $\phi$ such that the above limit exists
in $L^{2}(\Omega)$. To this end, we set\vadjust{\goodbreak} $I_{\varepsilon}(\phi
)=\int_{0}^{t}\dot{W}^{\varepsilon}(s,\phi_{s})\,ds$ and compute
$E (
I_{\varepsilon}(\phi)I_{\delta}(\phi) ) $ for \mbox{$\varepsilon
,\delta
>0$}. Denote
\[
V_{\varepsilon,\delta}^{2H} ( r ) =\frac{1}{4\varepsilon
\delta}%
 (  | r+\varepsilon-\delta | ^{2H}- |
r+\varepsilon
+\delta | ^{2H}- | r-\varepsilon-\delta |
^{2H}+ |
r-\varepsilon+\delta | ^{2H} ) .
\]
Using the fact that $Q ( x,y ) =Q ( y,x ) $, we have
\begin{eqnarray*}
&&E ( I_{\varepsilon}(\phi)I_{\delta}(\phi) )\\
 && \qquad =\frac{1}{
4\varepsilon\delta}\int_{0}^{t}\int_{0}^{\theta}Q(\phi_{\theta
},\phi
_{\eta})  [  | \theta-\eta+\varepsilon-\delta |
^{2H}- | \theta-\eta+\delta+\varepsilon | ^{2H} \\
&& \hspace*{98pt}\qquad  \quad  {}- | \theta-\eta-\varepsilon-\delta |
^{2H}+ |
\theta-\eta-\varepsilon+\delta | ^{2H} ] \,d\eta \,d\theta.
\end{eqnarray*}
Making the substitution $r={\theta}-\eta$ and using the notation $%
V_{\varepsilon,\delta}^{2H}$, we can write%
%
\begin{equation}
E ( I_{\varepsilon}(\phi)I_{\delta}(\phi) )
=\int_{0}^{t}\int_{0}^{\theta}Q(\phi_{\theta},\phi_{\theta
-r})V_{\varepsilon,\delta}^{2H} ( r )\, dr\,d\theta.
\label{e.3.square}
\end{equation}

We need the following two technical lemmas.

\begin{lemma}
\label{doubleLmt2}For any bounded function $\psi\dvtx [0,T]\rightarrow
\mathbb{R}
$, we have%
%
\begin{equation}
  \biggl|\int_{0}^{t}\psi ( s ) \int
_{0}^{s}V_{\varepsilon
,\delta}^{2H} ( r ) \,dr\,ds-2H\int_{0}^{t}\psi (
s )
s^{2H-1}\,ds \biggr|\leq4 \|\psi \|_{\infty}(\varepsilon+\delta)^{2H}.
\label{dctrl2}\hspace*{-35pt}
\end{equation}
\end{lemma}

\begin{pf}
$\!\!\!$Let $g ( s ) :=\int_{0}^{s} \vert r \vert
^{2H}\,dr$ and $%
f_{\varepsilon,\delta} ( t ) :=\int_{0}^{t}\psi (
s )
\int_{0}^{s}V_{\varepsilon,\delta}^{2H} ( r ) \,dr\,ds $.
Note that $%
g^{\prime\prime}$ exists everywhere except at $0$ and $g^{\prime
\prime
} ( r ) =2H  \operatorname{sign}   ( r )  \vert
r \vert
^{2H-1}$ for $r\neq0$. Then
\begin{eqnarray*}
f_{\varepsilon,\delta} ( t ) &=&\frac{1}{4\varepsilon
\delta} \int_{0}^{t}\psi ( s )   [ g ( s+\varepsilon
-\delta
 ) -g ( s+\varepsilon+\delta )\\
 &&\hphantom{\frac{1}{4\varepsilon
\delta} \int_{0}^{t}\psi ( s )   [}{} -g (
s-\varepsilon
-\delta ) +g ( s-\varepsilon+\delta )  ]\,ds \\
&=&\frac{1}{4}\int_{-1}^{1}\int_{-1}^{1}\int_{0}^{t}\psi (
s )
g^{\prime\prime} ( s+\eta\varepsilon-\xi\delta )
\,ds\,d\xi \,d\eta
\\
&=&\frac{1}{4}\int_{-1}^{1}\int_{-1}^{1}\int_{0}^{t}\psi (
s )
g^{\prime\prime} ( s-\Delta )\, ds\,d\xi \,d\eta,
\end{eqnarray*}
where $\Delta=\xi\delta-\eta\varepsilon$.

 \textit{Case} (i): If $\Delta\leq0$, we have
%
\begin{eqnarray}\label{e.3.4}
&& \biggl| \int_{0}^{t}\psi ( s )  \bigl( g^{\prime\prime
} (
s-\Delta ) -2Hs^{2H-1} \bigr)\,ds \biggr| \nonumber\\
&& \qquad \leq2H \| \psi \| _{\infty}\int
_{0}^{t} \bigl(
s^{2H-1}- ( s-\Delta ) ^{2H-1} \bigr)\,ds \\
&& \qquad = \| \psi \| _{\infty}  [
t^{2H}-(t-\Delta
)^{2H}+(-\Delta)^{2H}  ] \leq2 \| \psi \|
_{\infty} | \Delta | ^{2H}.\nonumber
\end{eqnarray}

\textit{Case} (ii): If $\Delta>0$, we assume that $\Delta<t$ (the
case $\Delta\geq t$ follows easily). Then
\begin{eqnarray*}
\int_{0}^{t}\psi ( s ) g^{\prime\prime} ( s-\Delta
 )\,ds&=&-2H\int_{0}^{\Delta}\psi ( s )  ( \Delta-s )
^{2H-1}\,ds\\[-3pt]
&&{}+2H\int_{\Delta}^{t}\psi ( s )  ( s-\Delta
 )
^{2H-1}\,ds .
\end{eqnarray*}
Therefore,
%
\begin{equation}
\biggl|\int_{0}^{t}\psi ( s )  \bigl( g^{\prime
\prime} (
s-\Delta ) -2Hs^{2H-1} \bigr)\,ds \biggr|\leq F_{\Delta
}^{1}+F_{\Delta}^{2} , \label{e.3.5}
\end{equation}
where
%
\begin{equation}
F_{\Delta}^{1}:=2H\int_{0}^{\Delta}\psi ( s )  [
 (
\Delta-s ) ^{2H-1}+s^{2H-1} ]\,ds\leq2\Vert\psi\Vert
_{\infty
}|\Delta|^{2H} \label{e.3.6}
\end{equation}
and
%
\begin{eqnarray} \label{e.3.7}
F_{\Delta}^{2} &:=&2H\int_{\Delta}^{t}\psi ( s )  [
 (
s-\Delta ) ^{2H-1}-s^{2H-1} ]\,ds \nonumber
\\[-9.5pt]
\\[-9.5pt]
&\hspace*{3pt}\leq&
2H\Vert\psi\Vert_{\infty}\int_{\Delta}^{t} [  (
s-\Delta
 ) ^{2H-1}-s^{2H-1} ]\,ds\leq2\Vert\psi\Vert_{\infty
}|\Delta
|^{2H} .
\nonumber
\end{eqnarray}
Then (\ref{dctrl2}) follows from (\ref{e.3.4})--(\ref{e.3.7}).\vspace*{-3pt}
\end{pf}

\begin{lemma}
\label{doubleLmt3} Let $\psi\in C( [ 0,T ] ^{2})$ with
$\psi
 ( 0,s ) =0,$ and $\psi ( \cdot,s ) \in
C^{\alpha}(%
 [ 0,T ] )$ for any $s\in [ 0,T ] $. Assume
$\alpha+2H>1$
and $\sup_{s\in\lbrack0,T]} \| \psi ( \cdot
,s )
 \| _{\alpha}<\infty$. Then for any $1-2H<{\gamma
}<{\alpha}$
and $t\leq T$ we have
%
\begin{eqnarray}
\label{e1}
&& \biggl| \int_{0}^{t}\int_{0}^{s}\psi ( r,s )  [
V_{\varepsilon,\delta}^{2H} ( r ) -2H ( 2H-1 )
r^{2H-2}%
 ] \,dr\,ds \biggr| \nonumber
 \\[-9.5pt]
 \\[-9.5pt]
&& \qquad \leq C\sup_{s\in\lbrack0,T]}\| \psi ( \cdot
,s )
 \| _{\alpha}({\varepsilon}+{\delta})^{2H+{\gamma
}-1} ,
\nonumber
\end{eqnarray}
where the constant $C$ depends on $H$, ${\gamma}$, $\alpha$ and $T$, but
it is independent of~${\delta},{\varepsilon}$ and~$\psi$.\vspace*{-3pt}
\end{lemma}

\begin{pf}
Along the proof, we denote by $C$ a generic constant which depends on
$H$, ${%
\gamma}$, $\alpha$ and $T$. Set $h ( r ) := \vert
r \vert^{2H}$. Then $h^{\prime}(r)$ exists everywhere except at $0$
and $h^{\prime} ( r ) =2H \operatorname{sign} ( r )
 \vert r \vert^{2H-1}$ if $r\neq0$. Using (\ref
{fracint1}) and %
(\ref{e.2.6}),  we have
\begin{eqnarray*}
 f_{\varepsilon,\delta} ( t ) &:=&\int_{0}^{t}\int
_{0}^{s}\psi
 ( r,s ) V_{\varepsilon\delta}^{2H} ( r ) \,dr\,ds
\\[-3pt]
&\hspace*{2.56pt}=&\frac{1}{4\varepsilon}\int_{-1}^{1}\int_{0}^{t}\int_{0}^{s}\psi
 (
r,s ) \,\frac{\partial}{\partial r} [ h ( r+\varepsilon
-\xi
\delta ) -h ( r-\varepsilon-\xi\delta )  ]
\,dr\,ds\,d\xi
\\[-3pt]
&\hspace*{2.56pt}=& ( -1 ) ^{\alpha^{\prime}}\frac{1}{4\varepsilon}%
\int_{-1}^{1}\int_{0}^{t}\int_{0}^{s}D_{0+}^{\alpha^{\prime}}\psi
 (
r,s ) \\[-3pt]
&&\hphantom{( -1 ) ^{\alpha^{\prime}}\frac{1}{4\varepsilon}%
\int_{-1}^{1}\int_{0}^{t}\int_{0}^{s}}{}\times D_{s-}^{1-\alpha^{\prime}}[ h ( r+\varepsilon
-\xi
\delta ) -h ( r-\varepsilon-\xi\delta )  ]
\,dr\,ds\,d\xi
\\[-3pt]
&\hspace*{2.56pt}=& ( -1 ) ^{\alpha^{\prime}}\frac{1}{4}\int_{-1}^{1}%
\int_{-1}^{1}\int_{0}^{t}\int_{0}^{s}D_{0+}^{\alpha^{\prime}}\psi
 (
r,s ) D_{s-}^{1-\alpha^{\prime}}h^{\prime} ( r+\eta
\varepsilon
-\xi\delta ) \,dr\,ds\,d\xi \,d\eta,
\end{eqnarray*}
where $\gamma<\alpha^{\prime}<\alpha$. On the other hand, we also have
\begin{eqnarray*}
&&2H(2H-1)\int_{0}^{t}\int_{0}^{s}\psi(r,s)r^{2H-2}\,dr\\
&& \qquad =(-1)^{{\alpha
}^{\prime
}}\int_{0}^{t}\int_{0}^{s}D_{0+}^{\alpha^{\prime}}\psi (
r,s )
D_{s-}^{1-\alpha^{\prime}}h^{\prime} ( r ) \,dr\,ds .
\end{eqnarray*}
Thus,
\begin{eqnarray*}
 I_{{\varepsilon},{\delta}}&:=& \biggl|\int_{0}^{t}\int
_{0}^{s}\psi
 ( r,s )  [ V_{\varepsilon,\delta}^{2H} (
r )
-2H ( 2H-1 ) r^{2H-2} ] \,dr\,ds \biggr| \\
&\hspace*{3pt}\leq&\frac{1}{4}\int_{-1}^{1}\int_{-1}^{1}\int_{0}^{t}\int
_{0}^{s} %
\vert D_{0+}^{\alpha^{\prime}}\psi ( r,s )  \vert\\
&&\hphantom{\frac{1}{4}\int_{-1}^{1}\int_{-1}^{1}\int_{0}^{t}\int
_{0}^{s}}{}\times
 \vert D_{s-}^{1-\alpha^{\prime}}h^{\prime} ( r+\eta
\varepsilon
-\xi\delta ) -D_{s-}^{1-\alpha^{\prime}}h^{\prime}(r)
\vert
\,dr\,ds\,d\xi \,d\eta.
\end{eqnarray*}
Denote $\Delta=\xi\delta-\eta\varepsilon$ and
%
\begin{equation}
 \qquad f_{\Delta} ( t ) :=\int_{0}^{t}\int_{0}^{s} \vert
D_{0+}^{\alpha^{\prime}}\psi ( r,s )  \vert
\vert
 [ D_{s-}^{1-\alpha^{\prime}}h^{\prime} ( r-\Delta )
-D_{s-}^{1-\alpha^{\prime}}h^{\prime} ( r )  ]
 \vert
\,dr\,ds. \label{fDelta}
\end{equation}
Then we may write
%
\begin{equation}
I_{{\varepsilon},{\delta}}\leq\frac{1}{4}\int_{-1}^{1}\int
_{-1}^{1}f_{{%
\Delta}}(t)\,d\xi \,d\eta. \label{Iepdel}
\end{equation}
Hence, in order to prove (\ref{e1}) it suffices to prove%
%
\begin{equation}
f_{\Delta} ( t ) \leq C\sup_{s\in\lbrack0,T]} \|\psi ( \cdot,s )  \|
_{\alpha}|{%
\Delta}|^{2H+{\gamma}-1} . \label{fDelta-}
\end{equation}
By (\ref{frDl}), we have
%
\begin{eqnarray}\label{frDphi}\qquad
 \vert D_{0+}^{\alpha^{\prime}}\psi ( r,s )
\vert &=&
\frac{1}{\Gamma ( 1-\alpha^{\prime} ) } \biggl|\frac
{\psi
 ( r,s ) }{r^{\alpha^{\prime}}}+\alpha^{\prime}\int
_{0}^{r}%
\frac{\psi ( r,s ) -\psi ( u,s ) }{ (
r-u )
^{\alpha^{\prime}+1}}\,du \biggr| \nonumber
\\[-8pt]
\\[-8pt]
&\leq&C\sup_{s\in\lbrack0,T]} \|\psi (
\cdot
,s )  \|_{\alpha}.
\nonumber
\end{eqnarray}
Therefore,%
%
\begin{equation}
f_{\Delta} ( t ) \leq C\sup_{s\in\lbrack0,T]} \|\psi ( \cdot,s )  \|
_{\alpha
} ( F_{\Delta}^{1}+F_{\Delta}^{2} ) , \label{a1}
\end{equation}
where
\begin{eqnarray*}
F_{\Delta}^{1} &=&\int_{0}^{t}\int_{0}^{s}\frac{ \vert
h^{\prime
} ( r-\Delta ) -h^{\prime} ( r )  \vert
}{ (
s-r ) ^{1-\alpha^{\prime}}}\,dr\,ds, \\
F_{\Delta}^{2} &=&\int_{0}^{t}\int_{0}^{s}\int_{r}^{s}\frac{
\vert
h^{\prime} ( r-\Delta ) -h^{\prime} ( u-\Delta
 )
-h^{\prime}(r)+h^{\prime}(u) \vert}{ ( u-r )
^{2-\alpha
^{\prime}}}\,du\,dr\,ds.
\end{eqnarray*}

As in the proof of Lemma \ref{doubleLmt2}, we consider the two cases
separately: ${\Delta}\leq0$ and ${\Delta}>0$.

 \textit{Case} (i): If $\Delta\leq0$, we can write
\begin{eqnarray*}
 \biggl|\frac{h^{\prime} ( r-\Delta ) -h^{\prime
}(r)}{ (
s-r ) ^{1-\alpha^{\prime}}} \biggr|&\leq&C(s-r)^{{\alpha}
^{\prime}-1}|{\Delta}|\int_{0}^{1}(r-\xi{\Delta})^{2H-2}\,d\xi  \\
&\leq&C(s-r)^{{\alpha}^{\prime}-1}r^{-\gamma}|{\Delta}|^{2H+\gamma
-1} ,
\end{eqnarray*}
which implies%
%
\begin{equation}
F_{\Delta}^{1}\leq C|{\Delta}|^{2H+\gamma-1}. \label{a2}
\end{equation}
For $0<r<u$, we have
\begin{eqnarray*}
&&   \vert h^{\prime} ( r-\Delta ) -h^{\prime} (
u-\Delta
 ) -h^{\prime}(r)+h^{\prime}(u) \vert\\
&& \qquad =C|{\Delta}|\int_{0}^{1}\int_{0}^{1} \bigl( r-\xi\Delta+\theta
 (
u-r )  \bigr) ^{2H-3}\,d\theta \,d\xi (r-u)  \\
&& \qquad \leq Cr^{2H-1-{\beta}_{1}-{\beta}_{2}}(u-r)^{{\beta}_{1}}|{\Delta
}|^{{%
\beta}_{2}}
\end{eqnarray*}
for any ${\beta}_{1},{\beta}_{2}>0$ such that ${\beta}_{1}+{\beta}%
_{2}<2H $. If ${\alpha}^{\prime}+{\beta}_{1}>1$, we obtain
\begin{eqnarray*}
&&\int_{r}^{s}\frac{|h^{\prime}(r-{\Delta})-h^{\prime}(u-{\Delta}%
)-h^{\prime}(r)+h(u)|}{(u-r)^{2-{\alpha}^{\prime}}}\,du  \\
&& \qquad \leq Cr^{2H-1-{\beta}_{1}-{\beta}_{2}}(s-r)^{{\alpha}^{\prime
}+{\beta}%
_{1}-1}|{\Delta}|^{{\beta}_{2}} ,
\end{eqnarray*}
which implies, taking $\beta_{2}=2H+{\gamma}-1$,%
%
\begin{equation}
F_{\Delta}^{2}\leq C|{\Delta}|^{2H+\gamma-1}. \label{a3}
\end{equation}
Substituting (\ref{a2}) and (\ref{a3}) into (\ref{a1}), we get (\ref
{fDelta-}%
).

 \textit{Case} (ii): Now let $\Delta>0$. We assume
that $\Delta<t$ (the case $t\leq\Delta$ is simpler and omitted). Let us
first consider the term $F_{\Delta}^{1}$. Define the sets
\begin{eqnarray*}
D_{11}&=&\{0<r<s<{\Delta}\} ,\qquad D_{12}=\{0<r<{\Delta} <s<t\}
,\\
D_{13}&=&\{{\Delta}<r<s <t\} .
\end{eqnarray*}
Then
\[
F_{\Delta}^{1}=F_{\Delta}^{11}+F_{\Delta}^{12}+F_{\Delta}^{13} ,
\]
where
\[
F_{\Delta}^{1i}=\int_{D_{1i}}\frac{ \vert h^{\prime} (
r-\Delta
 ) -h^{\prime} ( r )  \vert}{ ( s-r )
^{1-\alpha^{\prime}}}\,dr\,ds ,\qquad i=1,2,3 .\vadjust{\goodbreak}
\]
It is easy to see that
%
\begin{equation}
F_{\Delta}^{11}\leq C\int_{0}^{\Delta}\int_{0}^{s} [  (
\Delta
-r ) ^{2H-1}+r^{2H-1} ]  ( s-r ) ^{\alpha
^{\prime
}-1}\,dr\,ds\leq C\Delta^{2H+\alpha^{\prime}} \label{FDelta11}\hspace*{-35pt}
\end{equation}
and
%
\begin{equation}
F_{\Delta}^{12}\leq C\int_{\Delta}^{t}\int_{0}^{\Delta} [
 (
\Delta-r ) ^{2H-1}+r^{2H-1} ]  ( s-r ) ^{\alpha
^{\prime
}-1}\,dr\,ds\leq C\Delta^{2H} . \label{FDelta12}\hspace*{-35pt}
\end{equation}
As for $F_{\Delta}^{13}$, we have
\[
F_{\Delta}^{13}=\int_{\Delta}^{t}\int_{\Delta}^{s}\frac{
\vert
h^{\prime} ( r-\Delta ) -h^{\prime}(r) \vert
}{ (
s-r ) ^{1-\alpha}}\,dr\,ds=\int_{0}^{t-\Delta}\int_{0}^{u}\frac{%
 \vert h^{\prime} ( v ) -h^{\prime}(v+{\Delta
}) \vert}{%
 ( u-v ) ^{1-\alpha^{\prime}}}\,dv\,du .
\]
Using the estimate
\[
|h^{\prime}(v)-h^{\prime}(v+{\Delta})|\leq Cv^{2H-{\beta}-1}{\Delta
}^{{%
\beta}}
\]
for all $0<{\beta}<2H$, we obtain
%
\begin{equation}
F_{\Delta}^{13}\leq C {\Delta}^{{\beta}} . \label{FDelta13}
\end{equation}
Thus, (\ref{FDelta11})--(\ref{FDelta13}) yield
%
\begin{equation}
F_{\Delta}^{1}\leq C{\Delta}^{{\beta}}  \qquad\mbox{for
all } 0<{\beta}<2H . \label{FDelta1}
\end{equation}
Now we study the second term $F_{\Delta}^{2}$. Denote
\begin{eqnarray*}
  D_{21}&=&\{0<r<u<s<{\Delta<t}\} ,\qquad D_{22}=\{0<r<u<{\Delta}
<s<t\} ,
\\
  D_{23}&=&\{0<r<{\Delta} <u<s <t\} ,\qquad D_{24}=\{0<{\Delta}%
<r<u<s<t\} .
\end{eqnarray*}
Then
\[
F_{\Delta}^{2}= F_{\Delta}^{21}+F_{\Delta}^{22}+F_{\Delta
}^{23}+F_{\Delta}^{24} ,
\]
where for $i=1,2,3,4$,
\[
F_{\Delta}^{2i}=\int_{D_{2i}}\frac{ \vert h^{\prime} (
r-\Delta
 ) -h^{\prime} ( u-\Delta ) -h^{\prime
}(r)+h^{\prime
}(u) \vert}{ ( u-r ) ^{2-\alpha^{\prime}}}\,du\,dr\,ds .
\]
Consider first the term $F_{\Delta}^{21}$. We can write
\begin{eqnarray*}
 \frac{1}{2H} \vert h^{\prime} ( r-\Delta )
-h^{\prime
} ( u-\Delta )  \vert&=& \vert ( \Delta
-u )
^{2H-1}- ( \Delta-r ) ^{2H-1} \vert\\
&\leq&C ( u-r ) \int_{0}^{1} \bigl( \Delta-u+\theta (
u-r )  \bigr) ^{2H-2}\,d\theta\\
&\leq&C ( u-r ) ^{1-\beta} ( \Delta-u )
^{2H+\beta-2},
\end{eqnarray*}
where $1-2H<{\beta}<\alpha^{\prime}$. Similarly, we have
\[
 \vert h^{\prime}(r)-h^{\prime} ( u )  \vert
\leq
Cr^{2H+\beta-2}(u-r)^{1-{\beta}} .
\]
As a consequence,
%
\begin{eqnarray} \label{FDelta21}
F_{\Delta}^{21} &\leq& C\int_{0}^{\Delta}\int_{0}^{s}\int
_{r}^{s} (
u-r ) ^{\alpha^{\prime}-\beta-1} ( \Delta-u )
^{2H+\beta
-2}\,du\,dr\,ds \nonumber\\
&\leq& C\int_{0}^{\Delta}\int_{0}^{\Delta}\int_{r}^{\Delta} (
u-r ) ^{\alpha^{\prime}-\beta-1} ( \Delta-u )
^{2H+\beta
-2}\,du\,dr\,ds \\
&\leq& C\Delta^{2H+\alpha^{\prime}} .\nonumber
\end{eqnarray}
In a similar way we can prove that
%
\begin{eqnarray}\label{FDelta22}
  F_{\Delta}^{22}&\leq& C \int_{\Delta}^{t}\int_{0}^{\Delta}\int
_{r}^{\Delta
} ( u-r ) ^{\alpha^{\prime}-\beta-1} ( \Delta
-u )
^{2H+\beta-2}\,du\,dr\,ds\nonumber
\\[-8pt]
\\[-8pt]&\leq& C\Delta^{2H+\alpha^{\prime}-1}.
\nonumber
\end{eqnarray}
For $F_{\Delta}^{23}$, notice that when $r<\Delta<u,$%
\begin{eqnarray*}
 &&\vert h^{\prime} ( r-\Delta ) -h^{\prime} (
u-\Delta
 ) -h^{\prime}(r)+h^{\prime}(u) \vert\\
 && \qquad = ( \Delta
-r )
^{2H-1}+ ( u-\Delta ) ^{2H-1} +r^{2H-1}+u^{2H-1}
\end{eqnarray*}
and
\begin{eqnarray*}
 ( u-r ) ^{\alpha^{\prime}-2} &=& ( u-\Delta+\Delta
-r ) ^{\alpha^{\prime}-2} \\
&\leq& ( u-\Delta ) ^{-\beta} ( \Delta-r )
^{\alpha
^{\prime}+\beta-2}\wedge ( u-\Delta ) ^{-{\beta
}-2H+1} (
\Delta-r ) ^{2H+{\alpha}^{\prime}+{\beta}-3} ,
\end{eqnarray*}
where we can take any $\beta\in ( 0,1 ) $ satisfying
$2H+\beta
+\alpha^{\prime}>2.$ Then,
\begin{eqnarray*}
 F_{\Delta}^{23}&\leq& C\int_{D_{23}} [  ( \Delta-r )
^{2H-1}+ ( u-\Delta ) ^{2H-1}+r^{2H-1}+u^{2H-1} ]\\
&&\hphantom{C\int_{D_{23}}}{}\times
 (
u-r ) ^{\alpha^{\prime}-2}\,du\,dr\,ds \\
&\leq&C\int_{D_{23}} [  ( \Delta-r ) ^{2H+\alpha
^{\prime
}+\beta-3}( u-\Delta ) ^{-\beta}\\
&&\hphantom{C\int_{D_{23}} [}{} +r^{2H-1} (
u-\Delta
 ) ^{-\beta} ( \Delta-r ) ^{\alpha^{\prime}+\beta
-2}%
 ] \,du\,dr\,ds \\
&\leq&C\Delta^{2H+\alpha^{\prime}+\beta-2} .
\end{eqnarray*}
Taking $\beta=1+\gamma-\alpha^{\prime},$ we obtain
%
\begin{equation}
 \vert F_{\Delta}^{23} \vert\leq C\Delta^{2H+\alpha
^{\prime
}-1}. \label{FDelta23}
\end{equation}

Finally we consider the last term $F_{\Delta}^{24}$. Making the
substitutions $x=r-\Delta,$ $y=u-\Delta$ we can write
\begin{eqnarray*}
F_{\Delta}^{24} &=&\int_{D_{24}}\frac{ \vert h^{\prime} (
r-\Delta
 ) -h^{\prime} ( u-\Delta ) -h^{\prime
}(r)+h^{\prime
}(u) \vert}{ ( u-r ) ^{2-\alpha^{\prime}}}\,du\,dr\,ds \\
&=&\int_{\Delta}^{t}\int_{0}^{s-\Delta}\int_{x}^{s-\Delta}\frac{%
 \vert h^{\prime} ( x ) -h^{\prime} ( y )
-h^{\prime
}(x+{\Delta})+h^{\prime}(y+{\Delta}) \vert}{ ( y-x )
^{2-\alpha^{\prime}}}\,dy\,dx\,ds.
\end{eqnarray*}
Note that for $0<x<y$ and ${\Delta}>0$,
\begin{eqnarray*}
&& \vert h^{\prime} ( x ) -h^{\prime} ( y )
-h^{\prime}(x+{\Delta})+h^{\prime}(y+{\Delta}) \vert\\[-2pt]
&& \qquad =x^{2H-1}-y^{2H-1}-(x+{\Delta})^{2H-1}+(y+{\Delta})^{2H-1} \\[-2pt]
&& \qquad =C\int_{0}^{1}\int_{0}^{1}\bigl(x+{\theta}(y-x)+\tilde{{\theta
}}{\Delta}%
\bigr)^{2H-3}\,d{\theta}\,d\tilde{{\theta}} \\[-2pt]
&& \qquad \leq Cx^{2H+{\beta}_{1}+{\beta}_{2}-3}(y-x)^{1-{\beta}_{1}}{\Delta
}^{1-%
{\beta}_{2}} ,
\end{eqnarray*}
where
\[
0<{\beta}_{1},{\beta}_{2}<1 ,\qquad 2H+{\beta}_{1}+{\beta}_{2}>2
 \quad
\mbox{and}\quad{\beta}_{1}<{\alpha}^{\prime} .
\]
Taking ${\beta}_{2}=2-2H-\gamma$ we get
%
\begin{equation}
F_{{\Delta}}^{24}\leq C {\Delta}^{2H+{\gamma}-1} . \label{FDelta24}
\end{equation}
From (\ref{FDelta21})--(\ref{FDelta24}), we see that
%
\begin{equation}
F_{{\Delta}}^{2}\leq C{\Delta}^{2H+{\gamma}-1} . \label{FDelta2}
\end{equation}
This completes the proof of the lemma.
\end{pf}

\begin{theorem}
\label{def smooth approx} Suppose that $\phi \in C^{\alpha} (
 [
0,T ]  ) $ with $\gamma\alpha>1-2H$ on $ [ 0,T
] $.
Then, the nonlinear stochastic integral $\int_{0}^{t}W(ds,\phi_{s})$ exists
and
%
\begin{eqnarray}\label{Var}
 &&E \biggl( \int_{0}^{t}W(ds,\phi_{s}) \biggr) ^{2}\nonumber\hspace*{-30pt}\\[-2pt]
  && \qquad = 2H\int
_{0}^{t}\theta
^{2H-1}Q(\phi_{\theta},\phi_{\theta})\,d\theta\hspace*{-30pt}
\\[-2pt]
&& \qquad  \quad {}   +2H ( 2H-1 ) \int_{0}^{t}\int_{0}^{\theta
}r^{2H-2} \bigl(
Q(\phi_{\theta},\phi_{\theta-r})-Q(\phi_{\theta},\phi_{\theta
}) \bigr) \,dr\,d\theta.
\nonumber\hspace*{-30pt}
\end{eqnarray}
Furthermore, for any $\frac{1-2H}{\gamma} <{\alpha}^{\prime
}<{\alpha}$,
we have
%
\begin{eqnarray}\label{def int convg}
&&\sup_{0\leq t\leq T}E \biggl(  \biggl|\int_{0}^{t}\dot
{W}^{\varepsilon
}(s,\phi_{s})\,ds-\int_{0}^{t}W(ds,\phi_{s}) \biggr|^{2} \biggr)
\nonumber
\\[-9pt]
\\[-9pt]
&& \qquad \leq C ( 1+\Vert\phi\Vert_{\infty} ) ^{M} (
1+ \|\phi \|_{\alpha}^{\gamma} )
{%
\varepsilon}^{2H+\gamma{\alpha}^{\prime}-1},
\nonumber
\end{eqnarray}
where the constant $C$ depends on $H$, $T$, $\gamma$, $\alpha$,
$\alpha
^{\prime}$ and the constants $C_{0}$ and~$C_{1}$ appearing in (\textup{Q1})
and~(\textup{Q2}).
\end{theorem}

\begin{pf}
We can write  (\ref{e.3.square})  as
%
\begin{eqnarray}\label{e.3.29}
E ( I_{\varepsilon}(\phi)I_{\delta}(\phi) )
&=&\int_{0}^{t}\int_{0}^{\theta} \bigl( Q(\phi_{\theta},\phi
_{\theta
-r})-Q(\phi_{\theta},\phi_{\theta}) \bigr) V_{\varepsilon,\delta
}^{2H} ( r ) \,dr\,d\theta \nonumber
\\[-9pt]
\\[-9pt]
&&{}+\int_{0}^{t}\int_{0}^{\theta}Q(\phi_{\theta},\phi_{\theta
})V_{\varepsilon,\delta}^{2H} ( r ) \,dr\,d\theta .
\nonumber\vadjust{\goodbreak}
\end{eqnarray}
Due to the local boundedness of $Q$ [see (Q1)] and applying Lemma \ref%
{doubleLmt2} to $\psi ( \theta ) =Q(\phi_{\theta},\phi
_{\theta
})$, we see that the second integral converges to
\[
\lim_{{\varepsilon},{\delta}\rightarrow0}\int_{0}^{t}\int
_{0}^{\theta
}Q(\phi_{\theta},\phi_{\theta})V_{\varepsilon,\delta}^{2H} (
r ) \,dr\,d\theta=2H\int_{0}^{t}Q(\phi_{\theta},\phi_{\theta
})\theta
^{2H-1}\,d\theta .
\]
On the other hand, using the local H\"{o}lder continuity of $Q$ [see
(Q2)] and applying Lemma \ref{doubleLmt3}, to $\psi ( r,\theta
 )
=Q(\phi_{\theta},\phi_{\theta-r})-Q(\phi_{\theta},\phi_{\theta})$,
we see that the first integral converges to
\begin{eqnarray*}
&&\lim_{{\varepsilon},{\delta}\rightarrow0}\int_{0}^{t}\int
_{0}^{\theta
} \bigl( Q(\phi_{\theta},\phi_{\theta-r})-Q(\phi_{\theta},\phi
_{\theta
}) \bigr) V_{\varepsilon,\delta}^{2H} ( r ) \,dr\,d\theta\\
&& \qquad =2H ( 2H-1 ) \int_{0}^{t}\int_{0}^{\theta} \bigl(
Q(\phi_{\theta
},\phi_{\theta-r})-Q(\phi_{\theta},\phi_{\theta}) \bigr)
r^{2H-2}\,dr\,d\theta .
\end{eqnarray*}
This implies that $ \{ I_{\varepsilon_{n}}(\phi) ,n\geq1
\} $
is a Cauchy sequence in $L^{2}(\Omega)$ for any sequence ${\varepsilon
}%
_{n}\downarrow0$. As a consequence, $\lim_{\varepsilon\rightarrow
0}I_{\varepsilon}(\phi)$ exists in $L^{2}(\Omega)$ and is denoted by
$%
I(\phi):=\int_{0}^{t}W(ds,\phi_{s}) $. Letting ${\varepsilon
},{\delta}%
\rightarrow0$ in (\ref{e.3.29}), we obtain (\ref{Var}).

From  (\ref{e.3.29}), Lemma \ref{doubleLmt2} and Lemma \ref
{doubleLmt3}%
, we have for any ${\alpha}^{\prime}<{\alpha}$,
%
\begin{eqnarray}    \label{e.3.30}
&& \vert E ( I_{\varepsilon}(\phi)I_{\delta}(\phi) )
-E(I^{2}(\phi)) \vert\nonumber
\\[-8pt]
\\[-8pt]
&& \qquad \leq C ( 1+\Vert\phi\Vert_{\infty
} ) ^{M}(1+\Vert\phi\Vert_{\alpha}^{\gamma})({\varepsilon}+{
\delta})^{2H+\gamma{\alpha}^{\prime}-1} .
\nonumber
\end{eqnarray}
In equation {(\ref{e.3.30})}, let ${\delta}\rightarrow0$ and notice
that $I_{\delta}(\phi)\rightarrow I(\phi)$ in $L^{2}({\Omega})$. Then
\[
 \vert E ( I_{\varepsilon}(\phi)I(\phi) )
-E(I^{2}(\phi
)) \vert\leq C ( 1+\Vert\phi\Vert_{\infty} )
^{M}(1+\Vert
\phi\Vert_{\alpha}^{\gamma}){\varepsilon}^{2H+\gamma{\alpha
}^{\prime
}-1} .
\]
On the other hand, if we let ${\varepsilon}={\delta}$ in %
 {(\ref{e.3.30})}, we obtain
\[
 \vert EI_{\varepsilon}^{2}(\phi)-E(I^{2}(\phi)) \vert
\leq
C ( 1+\Vert\phi\Vert_{\infty} ) ^{M}(1+\Vert\phi\Vert
_{\alpha}^{\gamma}){\varepsilon}^{2H+\gamma{\alpha}^{\prime}-1} .
\]
Thus, we have
\[
E \vert I_{\varepsilon}(\phi)-I(\phi) \vert^{2}= [
E (
I_{\varepsilon}^{2}(\phi) ) -E ( I^{2}(\phi) )
 ] -2%
 [ E ( I_{\varepsilon}(\phi)I(\phi) ) -E (
I^{2}(\phi
) )  ]  .
\]
Applying the triangular inequality, we obtain  {(\ref{def
int convg})}.
\end{pf}

The following proposition can be proved in the same way as  {(\ref
{Var})}.

\begin{proposition}
Suppose $\phi,\psi\in C^{\alpha} ( [0,T] ) $ with
$\alpha
\gamma>1-2H$. Then
%
\begin{eqnarray}\label{Cov}
&& E \biggl( \int_{0}^{t}W(dr,\phi_{r})\int_{0}^{t}W(dr,\psi
_{r}) \biggr)\nonumber\\
&& \qquad =2H\int_{0}^{t}\theta^{2H-1}Q(\phi_{\theta},\psi_{\theta})\,d\theta
\nonumber
\\[-8pt]
\\[-8pt]
&& \qquad  \quad {}      +H ( 2H-1 ) \int_{0}^{t}\int_{0}^{\theta
}r^{2H-2} \bigl(
Q(\phi_{\theta},\psi_{\theta-r})-Q(\phi_{\theta},\psi_{\theta
}) \bigr) \,dr\,d\theta \nonumber\\
&& \qquad  \quad {}      +H ( 2H-1 ) \int_{0}^{t}\int_{0}^{\theta
}r^{2H-2} \bigl(
Q(\phi_{\theta-r},\psi_{\theta})-Q(\phi_{\theta},\psi_{\theta
}) \bigr)\,dr\,d\theta .\nonumber
\end{eqnarray}
\end{proposition}

The following proposition provides the H\"{o}lder continuity of the
indefinite integral.

\begin{proposition}
\label{nonlinear intgl lemma}Suppose $\phi \in C^{\alpha} (
[0,T] ) $ with $\alpha\gamma>1-2H$. Then for all $0\leq s<t\leq T$,
%
\begin{equation}
E\biggl ( \int_{0}^{t}W(dr,\phi_{r})-\int_{0}^{s}W(dr,\phi
_{r}) \biggr)
^{2} \leq C ( 1+\Vert\phi\Vert_{\infty} )
^{M}(t-s)^{2H} ,\hspace*{-35pt}
\label{e.3.35}
\end{equation}
where the constant $C$ depends on $H$, $T$, $\gamma$, $\alpha$ and the
constants $C_{0}$ and $C_{1}$ appearing in (\textup{Q1}) and (\textup{Q2}). As a
consequence, the process $X_{t}=\int_{0}^{t}W(dr,\phi_{r})$ is almost
surely $ ( H-\delta ) $-H\"{o}lder continuous for any
$\delta>0$.
\end{proposition}

\begin{pf}
We shall first show that%
%
\begin{equation}
  E\biggl ( \int_{0}^{t}\dot{W}^{\varepsilon}(r,\phi_{r})\,dr-\int
_{0}^{s}\dot{W}%
^{\varepsilon}(r,\phi_{r})\,dr \biggr) ^{2} \leq C ( 1+\Vert\phi
\Vert
_{\infty} ) ^{M}(t-s)^{2H} . \label{e.3.34}\hspace*{-35pt}
\end{equation}
We can write
\begin{eqnarray*}
 &&E \biggl( \int_{0}^{t}W^{\varepsilon}(dr,\phi
_{r})-\int_{0}^{s}W^{\varepsilon}(dr,\phi_{r}) \biggr) ^{2}\\
&& \qquad =E \biggl(
\int_{s}^{t}W^{\varepsilon}(dr,\phi_{r}) \biggr) ^{2} \\
&& \qquad =\frac{1}{4\varepsilon^{2}}\int_{s}^{t}\int_{s}^{t}E \bigl[ \bigl (
W(\theta+\varepsilon,\phi_{\theta})-W(\theta-\varepsilon,\phi
_{\theta
}) \bigr)\\
&&\hphantom{\frac{1}{4\varepsilon^{2}}\int_{s}^{t}\int_{s}^{t}E \bigl[ } \qquad  \quad {} \times\bigl ( W(\eta+\varepsilon,\phi_{\eta})-W(\eta
-\varepsilon
,\phi_{\eta}) \bigr)  \bigr]\,d\theta \,d\eta\\
&& \qquad =\frac{1}{8\varepsilon^{2}}\int_{s}^{t}\int_{s}^{t}Q(\phi
_{\theta},\phi
_{\eta}) \\
&& \hphantom{\frac{1}{8\varepsilon^{2}}\int_{s}^{t}\int_{s}^{t}}\qquad  \quad {}\times [  | \eta-\theta | ^{2H}- | \eta
-\theta
-2\varepsilon | ^{2H}- | \eta-\theta+2\varepsilon
| ^{2H}%
 ] \,d\theta \,d\eta\\
&& \qquad =\frac{1}{8\varepsilon^{2}}\int_{0}^{t-s}\int_{0}^{t-s}Q(\phi
_{s+\theta
},\phi_{s+\eta})\\
&&\hphantom{\frac{1}{8\varepsilon^{2}}\int_{0}^{t-s}\int_{0}^{t-s}} \qquad  \quad {}\times  [  | \eta-\theta | ^{2H}-
| \eta
-\theta-2\varepsilon | ^{2H}- | \eta-\theta+2\varepsilon
 | ^{2H} ] \,d\theta \,d\eta\\
&& \qquad =\frac{1}{4\varepsilon^{2}}\int_{0}^{t-s}\int_{0}^{\theta}Q(\phi
_{s+\theta},\phi_{s+\theta-r}) [ 2r^{2H}- | r+2\varepsilon
 | ^{2H}- | r-2\varepsilon | ^{2H} ] \,dr\,d\theta.
\end{eqnarray*}
The inequality (\ref{e.3.34}) follows from the assumption (Q1) and the
inequality~(\ref{DCT_ctrl}) obtained in the \hyperref[appm]{Appendix}. Finally, the
inequality~(\ref{e.3.35}) follows from~(\ref{e.3.34}), Proposition
\ref{def
smooth approx} and the Fatou's lemma.
\end{pf}

\section{Feynman--Kac integral}\label{sec4}

In this section, we show that the random field $u ( t,x ) $ given
by (\ref{e.1.3}) is well defined and study its H\"{o}lder continuity. Since
the Brownian motion $B_{t}$ has H\"{o}lder continuous trajectories of
order $%
\delta$ for any $\delta\in ( 0,\frac{1}{2} ) $, by Lemma
\ref%
{def smooth approx} the nonlinear stochastic integral $%
\int_{0}^{t}W(ds,B_{t-s}^{x})$ can be defined for any $H>\frac
{1}{2}-\frac{%
\gamma}{4} $. The following theorem shows that it is exponentially
integrable and hence $u ( t,x ) $ is well defined.

Set $ \| B \|_{\infty,T}=%
{\sup_{0\leq s\leq T}} \vert B_{s} \vert$ and $%
 \| B \|_{\delta,T}={%
\sup_{0\leq s<t\leq T}}\frac{|B_{t}-B_{s}|}{|t-s|^{\delta}}$ for
$\delta
\in ( 0,\frac{1}{2} ) .$

\begin{theorem}
\label{solu_integrability} Let $H>\frac{1}{2}-\frac{\gamma}{4} $
and let $%
u_{0}$ be bounded. For any $t\in\lbrack0,T]$ and $x\in\mathbb{R}^{d},$
the random variable $\int_{0}^{t}W(ds,B_{t-s}^{x})$ is exponentially
integrable and the random field $u ( t,x ) $ given by (\ref{e.1.3})
is in $L^{p}(\Omega)$ for any $p\geq1.$
\end{theorem}

\begin{pf}
Suppose first that $p=1$. By (\ref{e.3.35}) with $s=0$ and the Fernique's
theorem we have
\begin{eqnarray*}
E^{W} \vert u ( t,x )  \vert&\leq& \|
u_{0} \|_{\infty}E^{B}E^{W}\biggl[\exp
\int_{0}^{t}W(ds,B_{t-s}^{x})\biggr] \\
&\leq& \| u_{0} \|_{\infty
}E^{B}\bigl[e^{Ct^{2H}(1+ \| B \|
_{\infty
,T})^{M}}\bigr]<\infty.
\end{eqnarray*}

The $L^{p}$ integrability of $u ( t,x ) $ follows from Jensen's
inequality
%
\begin{eqnarray}\label{exp int of solu}
E^{W} | u ( t,x )  | ^{p} &\leq& \|
u_{0} \| _{\infty}E^{B}E^{W}\exp \biggl(
p\int_{0}^{t}W(dr,B_{t-r}^{x}) \biggr) \nonumber
\\[-8pt]
\\[-8pt]
&\leq& \| u_{0} \| _{\infty}E^{B}\bigl[\exp
 \bigl(Cp(1+ \| B \| _{\infty,T})^{M}T^{2H} \bigr)
\bigr]<\infty
 .
\nonumber
\end{eqnarray}
\upqed
\end{pf}

To show the H\"{o}lder continuity of $u ( \cdot,x ) $, we
need the
following lemma.

\begin{lemma}
\label{solu_hldr} Assume that $u_{0}$ is Lipschitz continuous. Then
for $%
0\leq s<t\leq T$ and for any $\alpha<2H-1+\frac{1}{2}\gamma$,
\[
E^{W} \biggl|
\int_{0}^{s}W(dr,B_{t-r}^{x})-\int_{0}^{s}W(dr,B_{s-r}^{x})
\biggr|
^{2}\leq C(1+ \| B \|
_{\infty
,T})^{M} \Vert B \Vert_{{\delta,T}}^{\gamma}(t-s)^{\alpha},
\]
where the constant $C$ depends on $H$, $T$, $\gamma$ and the constant $
C_{1} $ appearing in~(\textup{Q2}).
\end{lemma}

\begin{pf}
Suppose $\delta\in ( 0,\frac{1}{2} ) $. For $0\leq
u<v<s\leq T$,
denote%
\begin{eqnarray*}
\Delta Q ( s,t,u,v )
&:=&Q(B_{t-u}^{x},B_{t-v}^{x})-Q(B_{t-u}^{x},B_{t-u}^{x})\\
&&{}-Q(B_{t-u}^{x},B_{s-v}^{x})+Q(B_{t-u}^{x},B_{s-u}^{x}).
\end{eqnarray*}
Note that  (Q2)  implies
\[
 \vert\Delta Q ( s,t,u,v )  \vert\leq
2C_{1}(1+ \| B \|_{\infty
,T})^{M}\Vert B\Vert_{{\delta,T}}^{\gamma} ( t-s )
^{\gamma
\delta}
\]
and
\[
 \vert\Delta Q ( s,t,u,v )  \vert\leq
2C_{1}(1+ \| B \|_{\infty
,T})^{M} \Vert B \Vert_{{\delta,T}}^{\gamma} \vert
u-v \vert^{\gamma\delta},
\]
which imply that for any $\beta\in ( 0,1 ) ,$
\[
 \vert\Delta Q ( s,t,u,v )  \vert\leq
2C_{1}(1+ \| B \|_{\infty
,T})^{M} \Vert B \Vert_{{\delta,T}}^{\gamma} (
t-s )
^{\beta\gamma\delta} \vert u-v \vert^{ ( 1-\beta
 )
\gamma\delta}.
\]
Applying (\ref{Cov}) and using $Q ( x,y ) =Q (
y,x ), $ we
get
\begin{eqnarray*}
&&E^{W} \biggl|
\int_{0}^{s}W(dr,B_{t-r}^{x})-\int_{0}^{s}W(dr,B_{s-r}^{x})
\biggr|^{2}
\\
&& \qquad =2H ( 2H-1 ) \\
&& \qquad  \quad {}\times\int_{0}^{s}\int_{0}^{\theta}r^{2H-2}
[ \Delta
Q ( s,t,\theta,\theta-r ) +\Delta Q ( t,s,\theta
,\theta
-r )  ] \,dr\,d\theta\\
&& \qquad  \quad {}+2H\int_{0}^{s}\theta^{2H-1} \\
&&\hphantom{+2H\int_{0}^{s}} \qquad  \quad {}\times[ Q(B_{t-\theta}^{x},B_{t-\theta
}^{x})-2Q(B_{t-\theta}^{x},B_{s-\theta}^{x})+Q(B_{s-\theta
}^{x},B_{s-\theta}^{x}) ]\,d\theta\\
&& \qquad \leq C(1+ \| B \|_{\infty
,T})^{M}\Vert B\Vert_{{\delta,T}}^{\gamma}(t-s)^{\beta\gamma\delta}
\end{eqnarray*}
for any $\beta$ such that $ ( 1-\beta ) \gamma\delta>1-2H$,
that is, $\beta\gamma\delta<2H-1+\gamma\delta$. Taking~$\beta$ and
$\delta
$ such that $\beta\gamma\delta=\alpha$, we get the lemma.
\end{pf}

\begin{theorem}
\label{solu_Holder}Suppose $u_{0}$ is Lipschitz continuous and
bounded. Then
for each $x\in\mathbb{R}^{d}$, $u ( \cdot,x )  \in
C^{H_{1}} (  [ 0,T ]  ) $ for any $H_{1}\in
( 0,H-%
\frac{1}{2}+\frac{1}{4}\gamma ) $.
\end{theorem}

\begin{pf}
For $0\leq s<t\leq T$, from the Minkowski's inequality it follows that
%
\begin{eqnarray}\label{solu Holder 1}
&&E^{W} [  \vert u ( t,x ) -u ( s,x )
 \vert
^{p} ] \nonumber\\
&& \qquad \leq \bigl[ E^{B} \bigl( E^{W} \bigl|
u_{0}(B_{t}^{x})e^{\int_{0}^{t}W(dr,B_{t-r}^{x})}-u_{0}(B_{s}^{x})e^{%
\int_{0}^{s}W(dr,B_{s-r}^{x})} \bigr|^{p} \bigr) ^{\fraca
{1}{p}} \bigr]
^{p} \nonumber
\\[-8pt]
\\[-8pt]
&& \qquad \leq C \| u_{0} \|_{\infty
} \bigl[
E^{B} \bigl( E^{W} \bigl|
e^{\int_{0}^{t}W(dr,B_{t-r}^{x})}-e^{\int
_{0}^{s}W(dr,B_{s-r}^{x})} %
\bigr|^{p} \bigr) ^{{1}/{p}} \bigr] ^{p} \nonumber\\
&& \qquad  \quad {}+C \bigl[ E^{B} \bigl( E^{W} \bigl|
\bigl(u_{0}(B_{t}^{x})-u_{0}(B_{s}^{x})\bigr)e^{\int
_{0}^{s}W(dr,B_{s-r}^{x})} %
\bigr|^{p} \bigr) ^{{1}/{p}} \bigr] ^{p}.
\nonumber
\end{eqnarray}
Since $u_{0}$ is Lipschitz continuous, using (\ref{exp int of solu})
and H%
\"{o}lder's inequality, we have%
%
\begin{equation}
 \bigl[ E^{B} \bigl( E^{W} \bigl(  \vert
u_{0}(B_{t}^{x})-u_{0}(B_{s}^{x}) \vert
e^{\int_{0}^{s}W(dr,B_{s-r}^{x})} \bigr) ^{p} \bigr) ^{{1}/{p}} \bigr]
^{p}\leq C(t-s)^{{p}/{2}}. \label{solu Holder 0}\hspace*{-35pt}
\end{equation}
For the first term in (\ref{solu Holder 1}), using the formula that $%
|e^{a}-e^{b}|\leq(e^{a}+e^{b})|a-b|$ for $a,b\in\mathbb{R}$ and H\"
{o}%
lder's inequality we get
%
\begin{eqnarray}\label{solu Holder 2}
&&E^{W} \biggl[  \biggl|\exp\int_{0}^{t}W(dr,B_{t-r}^{x})-\exp
\int_{0}^{s}W(dr,B_{s-r}^{x}) \biggr|^{p} \biggr] \nonumber\\
&& \qquad \leq \biggl[ E^{W}\biggl( \exp\int_{0}^{t}W(dr,B_{t-r}^{x})+\exp
\int_{0}^{s}W(dr,B_{s-r}^{x})\biggr)^{2p} \biggr] ^{{1}/{2}} \\
&& \qquad  \quad {}  \times \biggl[ E^{W} \biggl|
\int_{0}^{t}W(dr,B_{t-r}^{x})-\int_{0}^{s}W(dr,B_{s-r}^{x})
\biggr|^{2p}%
 \biggr] ^{{1}/{2}}.\nonumber
\end{eqnarray}
Applying Lemma \ref{nonlinear intgl lemma} and Lemma \ref{solu_hldr}, we
obtain
%
\begin{eqnarray}
\label{solu_Holder3}
 &&E^{W} \biggl|
\int_{0}^{t}W(dr,B_{t-r}^{x})-\int_{0}^{s}W(dr,B_{s-r}^{x})
\biggr|^{2}\nonumber\\
&& \qquad \leq 2E^{W} \biggl|\int_{s}^{t}W(dr,B_{t-r}^{x}) \biggr|
^{2}
\nonumber
\\[-8pt]
\\[-8pt]
&& \qquad  \quad {}  +2E^{W} \biggl|
\int_{0}^{s}W(dr,B_{t-r}^{x})-\int_{0}^{s}W(dr,B_{s-r}^{x})
\biggr|^{2}
\nonumber\\
& & \qquad  \leq C(1+ \| B \|_{\infty
,T})^{M}\Vert B\Vert_{{\delta,T}}^{\gamma}(t-s)^{2H_{1}}.
\nonumber
\end{eqnarray}

Noting that conditional to $B$, $\int_{0}^{t}W(dr,B_{t-r}^{x})-%
\int_{0}^{s}W(dr,B_{s-r}^{x})$ is Gaussian, and using (\ref{solu
Holder 2}),
(\ref{solu_Holder3}) and (\ref{exp int of solu}) we get%
%
\begin{eqnarray}\label{solu Holder 4}
&& \biggl[ E^{B} \biggl( E^{W} \biggl|\exp\int
_{0}^{t}W(dr,B_{t-r}^{x})-\exp
\int_{0}^{s}W(dr,B_{s-r}^{x}) \biggr|^{p} \biggr) ^{{1}/{p}} \biggr]
^{p} \nonumber\\
&& \qquad \leq C\biggl [ E^{B} \biggl( E^{W} \biggl|
\int_{0}^{t}W(dr,B_{t-r}^{x})-\int_{0}^{s}W(dr,B_{s-r}^{x})
\biggr|
^{2} \biggr) ^{{1}/{2}} \biggr] ^{p} \\
&& \qquad \leq C(t-s)^{pH_{1}}.\nonumber
\end{eqnarray}
From (\ref{solu Holder 1}), (\ref{solu Holder 0}) and (\ref{solu
Holder 4}),
we can see that for any $p\geq1$,%
%
\begin{equation}
E^{W} [  \vert u ( t,x ) -u ( s,x )
 \vert
^{p} ] \leq C(t-s)^{pH_{1}}. \label{solu Holder 5}
\end{equation}
Now Kolmogorov's continuity criterion implies the theorem.
\end{pf}

\section{Validation of the Feynman--Kac formula}\label{sec5}

In the last section, we have proved that $u ( t,x ) $ given
by (\ref%
{e.1.3}) is well defined. In this section, we shall show that $u (
t,x ) $ is a weak solution to equation (\ref{she}).

To give the exact meaning about what we mean by a weak solution, we follow
the idea of \cite{H-N07} and \cite{hunuso}. First, we need a
definition of
the Stratonovich integral.

\begin{definition}
\label{def Str int} Given a random field $v=\{v(t,x)$, $t\geq0,x\in
\mathbb{%
R}^{d}\}$ such that $\int_{0}^{t}\int_{\mathbb
{R}^{d}}|v(s,x)|\,dx\,ds<\infty$
a.s. for all $t>0$, the Stratonovich integral
\[
\int_{0}^{t}\int_{\mathbb{R}^{d}}v(s,x)W(ds,x)\,dx
\]
is defined as the following limit in probability if it exists
\[
\lim_{\varepsilon\rightarrow0}\int_{0}^{t}\int_{\mathbb
{R}^{d}}v(s,x)\dot{W%
}^{\varepsilon}(s,x)\,ds\,dx,
\]
where $W^{\varepsilon}(t,x)$ is introduced in (\ref{W_approx}).
\end{definition}

The precise meaning of the weak solution to equation (\ref{she}) is given
below.

\begin{definition}
\label{def weak solu} A random field $u=\{u ( t,x ) ,t\geq
0,x\in
\mathbb{R}^{d}\}$ is a weak solution to equation (\ref{she}) if for
any $%
\varphi\in C_{0}^{\infty} ( \mathbb{R}^{d} ) $, we have
%
\begin{eqnarray}
\label{e.5.1}
\int_{\mathbb{R}^{d}}\bigl(u ( t,x ) -u_{0}(x)\bigr)\varphi(x)\,dx
&=&\int_{0}^{t}\int_{\mathbb{R}^{d}}u ( s,x ) \Delta
\varphi(x)\,dx\,ds \nonumber
\\[-8pt]
\\[-8pt]
&&{}+\int_{0}^{t}\int_{\mathbb{R}^{d}}u ( s,x ) \varphi(x)W(ds,x)\,dx
\nonumber\hspace*{-35pt}
\end{eqnarray}
almost surely, for all $t\geq0$, where the last term is a Stratonovich
stochastic integral in the sense of Definition \ref{def Str int}.
\end{definition}


The following theorem justifies the Feynman--Kac formula  {(\ref{e.1.3})}.

\begin{theorem}
\label{main_Th}Suppose $H>\frac{1}{2}-\frac{1}{4}\gamma$ and
$u_{0}$ is a
bounded measurable function. Let $u ( t,x ) $ be the random field
defined in (\ref{e.1.3}). Then for any $\varphi\in C_{0}^{\infty
} (
\mathbb{R}^{d} ) $, $u ( t,x ) \varphi(x)$ is Stratonovich
integrable and $u ( t,x ) $ is a weak solution to equation
(\ref%
{she}) in the sense of Definition \ref{def weak solu}.
\end{theorem}

\begin{pf}
We prove this theorem by a limit argument. We divide the proof into three
steps.

{\textit{Step}} 1.  Let $u^{\varepsilon}(t,x)$ be the unique
solution to the following equation:
%
\begin{equation}
\cases{\displaystyle
\frac{\partial u^{\varepsilon}}{\partial t}=\frac{1}{2}\Delta
u^{\varepsilon}+u^{\varepsilon}\,\frac{\partial W^{\varepsilon
}}{\partial t}%
(t,x),&\quad $t>0,x\in\mathbb{R}^{d} $,\vspace*{2pt}\cr\displaystyle
u^{\varepsilon}(0,x)=u_{0}(x) .%
}
\label{approxm}
\end{equation}
Since $W^{\varepsilon}(t,x)$ is differentiable, the classical Feynman--Kac
formula holds for the solution to this equation, that is,
\[
u^{\varepsilon} ( t,x ) :=E^{B}\bigl[u_{0}(B_{t}^{x})e^{\int
_{0}^{t}%
\dot{W}^{\varepsilon}(s,B_{t-s}^{x})\,ds}\bigr].
\]
The fact that $u^{\varepsilon} ( t,x ) $ is well defined follows
from (\ref{e.3.34}) and Fernique's theorem. In fact, we have (cf. the
argument in the proof of Lemma \ref{solu_integrability})
%
\begin{eqnarray}\label{uepsilon_int}
E^{W} \vert u^{\varepsilon} ( t,x )  \vert^{p}
&\leq& \| u_{0} \|_{\infty
}E^{B}E^{W}\exp
\biggl ( p\int_{0}^{t}\dot{W}^{\varepsilon}(r,B_{t-r}^{x})\,dr \biggr)
\nonumber
\\[-8pt]
\\[-8pt]
&\leq& \| u_{0} \|_{\infty
}E^{B}\bigl[\exp \bigl( Cp(1+ \| B \|
_{\infty,T})^{M}t^{2H} \bigr) \bigr]<\infty .
\nonumber
\end{eqnarray}

Introduce the following notations
\begin{eqnarray*}
g_{s,x}^{\varepsilon}(r,z)&:=&  \frac{1}{2\varepsilon}\mathbf{1}%
_{[s-\varepsilon, s+\varepsilon]}(r)\mathbf{1}_{[0,x]}(z), \\
g_{s,x}^{B} ( r,z ) &:=&  \mathbf{1}_{[0,s]}(r)\mathbf{1}%
_{[0,B_{s-r}^{x}]}(z), \\
g_{s,x}^{\varepsilon,B} ( r,z ) &:=&  \int_{0}^{s}\frac{1}{
\varepsilon}\mathbf{1}_{[\theta-\varepsilon, \theta+\varepsilon
]}(r)%
\mathbf{1}_{[0,B_{s-\theta}^{x}]}(z)\,d\theta .
\end{eqnarray*}
From the results of Section~\ref{sec3}, we see that $g_{s,x}^{\varepsilon}$, $%
g_{s,x}^{B}$, $g_{s,x}^{\varepsilon,B}\in\mathcal{H}$ ($\mathcal
{H}$ is
introduced in Section~\ref{sec2}), and we can write%
\begin{eqnarray*}
\dot{W}^{\varepsilon}(s,x)&=&W \biggl( \frac{1}{2\varepsilon}\mathbf
{1}%
_{[s-\varepsilon, s+\varepsilon]}(r)\mathbf{1}_{[0,x]}(z) \biggr)
=W (
g_{s,x}^{\varepsilon} ) ,
\\
\int_{0}^{s}W ( d\theta,B_{s-\theta}^{x} ) &=&W (
g_{s,x}^{B} ) ,\qquad\int_{0}^{s}\dot{W}^{\varepsilon} (
\theta
,B_{s-\theta}^{x} )\,d\theta=W ( g_{s,x}^{\varepsilon
,B} ) .
\end{eqnarray*}
Set
\[
\widetilde{u}^{\varepsilon}(s,x):=u^{\varepsilon} ( s,x )
-u ( s,x ) .
\]

{\textit{Step}} 2.  We prove the following claim:

$u^{\varepsilon}(s,x)\rightarrow u(s,x)$ in $\mathbb{D}^{1,2}$ as $%
\varepsilon\downarrow0$, uniformly on any compact subset of
${[0,T]\times
\mathbb{R}^{d}}$, that is, for any compact ${K}\subseteq{\mathbb{R}^{d}}$
%
\begin{equation}\qquad
\sup_{s\in\lbrack0,T],x\in K}E^{W} [ |\widetilde
{u}^{\varepsilon
}(s,x)|^{2}+ \Vert D\widetilde{u}^{\varepsilon}(s,x) \Vert
_{%
\mathcal{H}}^{2} ] \rightarrow0 \qquad \mbox{as }\varepsilon
\downarrow0.
\label{e.5.6}
\end{equation}

Since $u_{0}$ is bounded, without loss of generality, we may assume $%
u_{0}\equiv1$. Let $B^{1}$ and $B^{2}$ be two independent Brownian motions,
both independent of~$W$. Using the inequality $|e^{a}-e^{b}|\leq
(e^{a}+e^{b})|a-b|$, H\"{o}lder inequality and the fact that $%
W(g_{t,x}^{\varepsilon,B})$ and $W(g_{t,x}^{B})$ are Gaussian conditioning
to $B$, we have
\begin{eqnarray*}
 &&E^{W} \bigl( u^{\varepsilon} ( t,x ) -u ( t,x
)  \bigr)
^{2}\\
&& \qquad =E^{W}\bigl[E^{B}\bigl(e^{W(g_{t,x}^{\varepsilon
,B})}-e^{W(g_{t,x}^{B})}\bigr)\bigr]^{2} \\
&& \qquad \leq E^{B}E^{W} \bigl| e^{W(g_{t,x}^{\varepsilon
,B})}-e^{W(g_{t,x}^{B})} \bigr|^{2} \\
&& \qquad \leq E^{B} \bigl[ E^{W} \bigl( e^{W(g_{t,x}^{\varepsilon
,B})}+e^{W(g_{t,x}^{B})} \bigr) ^{4} \bigr] ^{{1}/{2}} [
E^{W} \vert W(g_{t,x}^{\varepsilon,B})-W(g_{t,x}^{B}) \vert
^{4}%
 ] ^{{1}/{2}} \\
&& \qquad \leq C \bigl[ E^{B}E^{W} \bigl( e^{4W(g_{t,x}^{\varepsilon
,B})}+e^{4W(g_{t,x}^{B})} \bigr)  \bigr] ^{{1}/{2}}E^{B}E^{W} \vert
W(g_{t,x}^{\varepsilon,B})-W(g_{t,x}^{B}) \vert^{2}.
\end{eqnarray*}
Note that (\ref{exp int of solu}) and (\ref{uepsilon_int}) imply%
%
\begin{equation}
E^{B}E^{W} \bigl( e^{pW(g_{t,x}^{\varepsilon
,B})}+e^{pW(g_{t,x}^{B})} \bigr)
<\infty \label{exp_int}
\end{equation}
for any $p\geq1.$ On the other hand, applying Theorem \ref{def smooth
approx}, we have
%
\begin{equation}
\sup_{0\leq t\leq T,x\in K}E^{B}E^{W} \vert
W(g_{t,x}^{\varepsilon
,B})-W(g_{t,x}^{B}) \vert^{2}\rightarrow0 \qquad \mbox{as
}\varepsilon
\downarrow0. \label{uepsilon_2}
\end{equation}
Then it follows that as $\varepsilon\downarrow0$
\[
\sup_{0\leq t\leq T,x\in K}E^{W} \vert\widetilde
{u}^{\varepsilon
} ( t,x )  \vert^{2}=\sup_{0\leq t\leq T,x\in
K}E^{W}\bigl(u^{\varepsilon} ( t,x ) -u ( t,x )
\bigr)^{2}\rightarrow0.
\]
For the Malliavin derivatives, we have%
\begin{eqnarray*}
Du^{\varepsilon}(s,x) &=&E^{B} [ \exp ( W (
g_{s,x}^{\varepsilon,B} )  ) g_{s,x}^{\varepsilon
,B} ] , \\
Du(s,x) &=&E^{B} [ \exp ( W ( g_{s,x}^{B} )
 )
g_{s,x}^{B} ] .
\end{eqnarray*}
Then
\begin{eqnarray*}
&&E^{W} \| Du^{\varepsilon}(s,x)-Du(s,x)
\|_{\mathcal{H}}^{2} \\
&& \qquad =E^{W} \bigl\| E^{B} \bigl[  \bigl( \exp (
W (
g_{s,x}^{\varepsilon,B} )  ) g_{s,x}^{\varepsilon,B}-\exp
 (
W ( g_{s,x}^{B} )  ) g_{s,x}^{B} \bigr)  \bigr]
 \bigr\|_{\mathcal{H}}^{2} \\
&& \qquad \leq2E^{W}E^{B} [ \exp ( 2W ( g_{s,x}^{\varepsilon
,B} )  )  \| g_{s,x}^{\varepsilon
,B}-g_{s,x}^{B} \|_{\mathcal{H}}^{2} ] \\
&& \qquad  \quad {}+2E^{W}E^{B} [  \vert\exp ( W (
g_{s,x}^{\varepsilon
,B} )  ) -\exp ( W ( g_{s,x}^{B} )  )
 \vert^{2} \| g_{s,x}^{B} \|_{%
\mathcal{H}}^{2} ] .
\end{eqnarray*}
Note that $\|g_{t,x}^{\varepsilon,B}-g_{t,x}^{B}\|_{\mathcal{H}%
}^{2}=E^{W} \vert W(g_{t,x}^{\varepsilon,B})-W(g_{t,x}^{B})
\vert
^{2}$. Then it follows again from (\ref{exp_int}) and (\ref{uepsilon_2})
that as $\varepsilon\downarrow0$%
\[
\sup_{0\leq t\leq T,x\in K}E^{W} \|
Du^{\varepsilon
}(s,x)-Du(s,x) \|_{\mathcal{H}}^{2}\rightarrow0.
\]
{\textit{Step}} 3. From equation (\ref{approxm}) and (\ref{e.5.6}%
), it follows that $\int_{0}^{t}\int_{\mathbb{R}^{d}}u^{\varepsilon
} (
s,x ) \varphi(x)\times \dot{W}^{\varepsilon}(s,x)\,ds\,dx$ converges in $L^{2}$
to some random variable as $\varepsilon\downarrow0.$ Hence, if
%
\begin{equation}
V_{\varepsilon}:=\int_{0}^{t}\int_{\mathbb{R}^{d}}\bigl(u^{\varepsilon
} (
s,x ) -u ( s,x ) \bigr)\varphi(x)\dot{W}^{\varepsilon
}(s,x)\,ds\,dx
\label{Vepsilon}
\end{equation}
converges to zero in $L^{2}$, then
\begin{eqnarray*}
&&\lim_{\varepsilon\rightarrow0}\int_{0}^{t}\int_{\mathbb
{R}^{d}}u (
s,x ) \varphi(x)\dot{W}^{\varepsilon}(s,x)\,ds\,dx\\
&& \qquad =\lim
_{\varepsilon
\rightarrow0}\int_{0}^{t}\int_{\mathbb{R}^{d}}u^{\varepsilon
}(s,x)\varphi
(x)\dot{W}^{\varepsilon}(s,x)\,ds\,dx,
\end{eqnarray*}
that is, $u ( s,x ) \varphi(x)$ is Stratonovich integrable
and $%
u ( s,x ) $ is a weak solution to equation (\ref{she}).
Thus, it
remains to show that $V_{\varepsilon}$ converges to zero in $L^{2}$.

  In order to show the convergence to zero of (\ref{Vepsilon}) in $L^{2}$,
first we write $\widetilde{u}^{\varepsilon
}(s,x)W(g_{s,x}^{\varepsilon})$
as the sum of a divergence integral and a trace term [see~(\ref{div_fmla})]
\[
\widetilde{u}^{\varepsilon}(s,x)W(g_{s,x}^{\varepsilon})=\delta(%
\widetilde{u}^{\varepsilon}(s,x)g_{s,x}^{\varepsilon})- \langle
D%
\widetilde{u}^{\varepsilon}(s,x),g_{s,x}^{\varepsilon} \rangle
_{%
\mathcal{H}}.
\]
Then we have
\begin{eqnarray*}
V_{\varepsilon} &=&\int_{0}^{t}\int_{\mathbb{R}^{d}}\widetilde{u}%
^{\varepsilon}(s,x)\varphi(x)W(g_{s,x}^{\varepsilon})\,ds\,dx  \\
&=&\int_{0}^{t}\int_{\mathbb{R}^{d}} \bigl( \delta(\widetilde{u}%
^{\varepsilon}(s,x)g_{s,x}^{\varepsilon})- \langle D\widetilde
{u}%
^{\varepsilon}(s,x),g_{s,x}^{\varepsilon} \rangle_{\mathcal{H}
} \bigr) \varphi(x)\,ds\,dx \\
&=&\delta(\psi^{\varepsilon})-\int_{0}^{t}\int_{\mathbb{R}%
^{d}} \langle D\widetilde{u}^{\varepsilon
}(s,x),g_{s,x}^{\varepsilon
} \rangle_{\mathcal{H}}\varphi(x)\,ds\,dx=:V_{\varepsilon
}^{1}-V_{\varepsilon}^{2},
\end{eqnarray*}
where
\[
\psi^{\varepsilon}(r,z)=\int_{0}^{t}\int_{\mathbb
{R}^{d}}\widetilde{u}%
^{\varepsilon}(s,x)g_{s,x}^{\varepsilon}(r,z)\varphi(x)\,ds\,dx.
\]
For the term $V_{\varepsilon}^{1}$, using the estimates on $L^{2}$
norm of
the Skorohod integral (see  (1.47)  in \cite{Nu06}), we
obtain%
%
\begin{equation}
E [  \vert V_{\varepsilon}^{1} \vert^{2} ] \leq
E [
 \|\psi^{\varepsilon} \|
_{\mathcal{%
H}}^{2} ] +E [  \| D\psi^{\varepsilon
} \|_{\mathcal{H}\otimes\mathcal
{H}}^{2} ] .
\label{div_term}
\end{equation}
Denoting $\operatorname{supp} ( \varphi ) $ the support of
$\varphi$,
we have
\begin{eqnarray*}
&&E [  \|\psi^{\varepsilon} \|
_{\mathcal{H}}^{2} ] \\
&& \qquad =E\int_{0}^{t}\int_{0}^{t}\int_{\mathbb{R}^{d}}\int_{\mathbb
{R}^{d}}%
\widetilde{u}^{\varepsilon}(s_{1},x_{1})\widetilde{u}^{\varepsilon
}(s_{2},x_{2}) \\
&& \qquad  \quad \hphantom{E\int_{0}^{t}\int_{0}^{t}\int_{\mathbb{R}^{d}}\int_{\mathbb
{R}^{d}}}{}\times \langle g_{s_{1},x_{1}}^{\varepsilon
},g_{s_{2},x_{2}}^{\varepsilon} \rangle_{\mathcal{H}}\varphi
(x_{1})\varphi(x_{2})\,ds_{1}\,ds_{2}\,dx_{1}\,dx_{2} \\
&& \qquad \leq M_1 \int_{0}^{t}%
\int_{0}^{t}\int_{\mathbb{R}^{d}}\int_{\mathbb{R}^{d}}
\langle
g_{s_{1},x_{1}}^{\varepsilon},g_{s_{2},x_{2}}^{\varepsilon}
\rangle_{%
\mathcal{H}}\varphi(x_{1})\varphi(x_{2})\,ds_{1}\,ds_{2}\,dx_{1}\,dx_{2} \\
&& \qquad =M_1 \int_{\mathbb{R}^{d}}\int_{%
\mathbb{R}^{d}} E^{W} [ W^{\varepsilon} ( t,x_{1} )
W^{\varepsilon} ( t,x_{2} )  ] \varphi(x_{1})\varphi
(x_{2})\,dx_{1}\,dx_{2},
\end{eqnarray*}
where $M_1:=\sup_{s\in\lbrack0,T],x\in\operatorname{supp} ( \varphi
 ) }E [
|\widetilde{u}^{\varepsilon}(s,x)|^{2} ] $. Note that
%
\begin{eqnarray}\label{cnvg_gepsilon}
&&\lim_{\varepsilon\rightarrow0}\int_{\mathbb{R}^{d}}\int_{\mathbb
{R}%
^{d}} E^{W} [ W^{\varepsilon} ( t,x_{1} )
W^{\varepsilon
} ( t,x_{2} )  ] \varphi(x_{1})\varphi(x_{2})\,dx_{1}\,dx_{2}
\nonumber\\
&& \qquad =\int_{\mathbb{R}^{d}}\int_{\mathbb{R}^{d}} E^{W} [ W (
t,x_{1} ) W ( t,x_{2} )  ] \varphi
(x_{1})\varphi
(x_{2})\,dx_{1}\,dx_{2} \\
&& \qquad =\int_{\mathbb{R}^{d}}\int_{\mathbb{R}^{d}} t^{2H}Q (
x_{1},x_{2} ) \varphi(x_{1})\varphi(x_{2})\,dx_{1}\,dx_{2}<\infty .\nonumber
\end{eqnarray}
Thus by  {(\ref{e.5.6})}, we get $E [  \|\psi
^{\varepsilon} \|_{\mathcal{H}}^{2} ]
\rightarrow
0$ as $\varepsilon\downarrow0$.\vadjust{\goodbreak}

On the other hand, setting $M_2:=\sup_{s\in\lbrack0,T],x\in\mathrm
{supp} ( \varphi ) }E%
 [ \| D\widetilde{u}^{\varepsilon}(s,x) \|_{\mathcal{%
H}}^{2} ] $, we have
\begin{eqnarray*}
&&E [  \| D\psi^{\varepsilon} \|
_{\mathcal{H}\otimes\mathcal{H}}^{2} ] \\
&& \qquad =E\int_{0}^{t}\int_{0}^{t}\int_{\mathbb{R}^{d}}\int_{\mathbb{R}%
^{d}} \langle D\widetilde{u}^{\varepsilon}(s_{1},x_{1})\otimes
g_{s_{1},x_{1}}^{\varepsilon},D\widetilde{u}^{\varepsilon
}(s_{2},x_{2})\otimes g_{s_{2},x_{2}}^{\varepsilon} \rangle
_{\mathcal{%
H}} \\
&& \hphantom{E\int_{0}^{t}\int_{0}^{t}\int_{\mathbb{R}^{d}}\int_{\mathbb{R}%
^{d}} }\qquad  \quad {}\times\varphi(x_{1})\varphi(x_{2})\,ds_{1}\,ds_{2}\,dx_{1}\,dx_{2} \\
&& \qquad =E\int_{0}^{t}\int_{0}^{t}\int_{\mathbb{R}^{d}}\int_{\mathbb{R}%
^{d}} \langle D\widetilde{u}^{\varepsilon
}(s_{1},x_{1}),D\widetilde{u}%
^{\varepsilon}(s_{2},x_{2}) \rangle_{\mathcal{H}}
\langle
g_{s_{1},x_{1}}^{\varepsilon},g_{s_{2},x_{2}}^{\varepsilon}
\rangle_{%
\mathcal{H}}\\
&& \hphantom{E\int_{0}^{t}\int_{0}^{t}\int_{\mathbb{R}^{d}}\int_{\mathbb{R}%
^{d}} }\qquad  \quad {}\times
\varphi(x_{1})\varphi(x_{2})\,ds_{1}\,ds_{2}\,dx_{1}\,dx_{2} \\
&& \qquad \leq M_2 \int_{0}^{t}\int_{0}^{t}\int_{\mathbb{R}^{d}}\int
_{\mathbb{R}%
^{d}}  \langle g_{s_{1},x_{1}}^{\varepsilon
},g_{s_{2},x_{2}}^{\varepsilon} \rangle_{\mathcal{H}}\varphi
(x_{1})\varphi(x_{2})\,ds_{1}\,ds_{2}\,dx_{1}\,dx_{2}.
\end{eqnarray*}
Then  {(\ref{e.5.6})} and (\ref{cnvg_gepsilon}) imply that
$E [
 \| D\psi^{\varepsilon} \|
_{%
\mathcal{H}\otimes\mathcal{H}}^{2} ] $ converges to zero as \mbox{$%
\varepsilon\downarrow0$}.

Finally, we deal with the trace term
%
\begin{eqnarray}\label{Trace}
   V_{\varepsilon}^{2} &\hspace*{-3pt}=&\int_{0}^{t}\int_{\mathbb
{R}^{d}} \bigl(  \langle
Du^{\varepsilon} ( s,x ) ,g_{s,x}^{\varepsilon}
\rangle_{%
\mathcal{H}}- \langle Du ( s,x ) ,g_{s,x}^{\varepsilon
} \rangle_{\mathcal{H}} \bigr) \varphi(x)\,ds\,dx \nonumber\hspace*{-35pt}
\\[-8pt]
\\[-8pt]
&=:& T_{1}^{\varepsilon}-T_{2}^{\varepsilon},
\nonumber\hspace*{-35pt}
\end{eqnarray}
where%
\begin{eqnarray*}
T_{1}^{\varepsilon} &=&\int_{0}^{t}\int_{\mathbb{R}^{d}}
\langle
Du^{\varepsilon} ( s,x ) ,g_{s,x}^{\varepsilon}
\rangle_{%
\mathcal{H}}\varphi(x)\,ds\,dx, \\
T_{2}^{\varepsilon} &=&\int_{0}^{t}\int_{\mathbb{R}^{d}}
\langle
Du ( s,x ) ,g_{s,x}^{\varepsilon} \rangle_{\mathcal
{H}%
}\varphi(x)\,ds\,dx.
\end{eqnarray*}
We will show that $T_{1}^{\varepsilon}$ and $T_{2}^{\varepsilon}$ converge
to the same random variable as $\varepsilon\downarrow0$.

We start with the term $T_{2}^{\varepsilon}$. Note that
\begin{eqnarray*}
 \langle g_{s,x}^{B},g_{s,x}^{\varepsilon} \rangle&=&
\biggl\langle
\mathbf{1}_{[0,s]}(r)\mathbf{1}_{[0,B_{s-r}^{x}]}(z),\frac
{1}{2\varepsilon}%
\mathbf{1}_{[s-\varepsilon, s+\varepsilon]}(r)\mathbf{1}%
_{[0,x]}(z) \biggr\rangle_{\mathcal{H}} \\
&=& \biggl\langle\mathbf{1}_{[0,s]}(r)Q ( B_{s-r}^{x},x )
,\frac{1}{%
2\varepsilon}\mathbf{1}_{[s-\varepsilon, s+\varepsilon]}(r)
\biggr\rangle
_{\mathcal{H}}.
\end{eqnarray*}
Since $Q ( B_{s-\cdot}^{x},x ) \in C^{{1}/{2}-\delta
} ( %
 [ 0,T ]  ) $ for any $0<\delta<\frac{1}{2}$,
noticing that $%
H>\frac{1}{2}-\frac{\gamma}{4}$ and applying Lemma \ref{DCT_Lm2} we obtain
%
\begin{eqnarray}\label{Trace2}
\lim_{\varepsilon\rightarrow0}T_{2}^{\varepsilon
}&=&E^{B}\int_{0}^{t}\int_{\mathbb{R}^{d}}u_{0} (
B_{s}^{x} ) \exp
 ( W ( g_{s,x}^{B} )  )  \langle
g_{s,x}^{B},g_{s,x}^{\varepsilon} \rangle_{\mathcal{H}}\varphi(x)\,ds\,dx
\nonumber\\
&=&E^{B}\int_{0}^{t}\int_{\mathbb{R}^{d}}u_{0} ( B_{s}^{x}
) \exp
 ( W ( g_{s,x}^{B} )  ) \varphi(x) \nonumber\\
&&\hphantom{E^{B}\int_{0}^{t}\int_{\mathbb{R}^{d}}}{}\times\biggl[Q (
x,x ) Hs^{2H-1}\\
&&\hphantom{E^{B}\int_{0}^{t}\int_{\mathbb{R}^{d}}{}\times\biggl[}{}
+H ( 2H-1 ) \int_{0}^{s} \bigl( Q (
B_{s-r}^{x},x )\nonumber\\
&&\hspace*{155pt}{}
-Q ( x,x )  \bigr) r^{2H-2}\,dr\biggr]\,ds\,dx.\nonumber
\end{eqnarray}
On the other hand, for the term $T_{1}^{\varepsilon}$, note that
\begin{eqnarray*}
 \langle g_{s,x}^{\varepsilon,B},g_{s,x}^{\varepsilon}
\rangle
&=& \biggl\langle\int_{0}^{2\varepsilon}\frac{1}{2\varepsilon
}\mathbf{1}%
_{[\theta-\varepsilon,\theta+\varepsilon]}(r)\mathbf
{1}_{[0,B_{s-\theta
}^{x}]}(z)\,d\theta,\frac{1}{2\varepsilon}\mathbf{1}_{[s-\varepsilon
, s+\varepsilon]}(r)\mathbf{1}_{[0,x]}(z) \biggr\rangle_{\mathcal
{H}} \\
&=& \biggl\langle\int_{0}^{2\varepsilon}\frac{1}{2\varepsilon
}\mathbf{1}%
_{[\theta-\varepsilon,\theta+\varepsilon]}(r)Q ( B_{s-\theta
}^{x},x )\,d\theta,\frac{1}{2\varepsilon}\mathbf
{1}_{[s-\varepsilon
, s+\varepsilon]}(r) \biggr\rangle_{\mathcal{H}}.
\end{eqnarray*}

\mbox{}

\noindent Applying Lemma \ref{DCT_Lm4}, we obtain
%
\begin{eqnarray}\label{Trace1}
\lim_{\varepsilon\rightarrow0}T_{1}^{\varepsilon
}&=&E^{B}\int_{0}^{t}\int_{\mathbb{R}^{d}}u_{0} (
B_{s}^{x} ) \exp
 ( W ( g_{s,x}^{\varepsilon,B} )  )
\langle
g_{s,x}^{\varepsilon,B},g_{s,x}^{\varepsilon} \rangle
_{\mathcal{H}%
}\varphi(x)\,ds\,dx  \nonumber\\
&=&E^{B}\int_{0}^{t}\int_{\mathbb{R}^{d}}u_{0} ( B_{s}^{x}
) \exp
 ( W ( g_{s,x}^{B} )  ) \varphi(x) \nonumber\\
&&\hphantom{E^{B}\int_{0}^{t}\int_{\mathbb{R}^{d}}}{}\times\biggl[Q (
x,x ) Hs^{2H-1}\\
&&\hphantom{E^{B}\int_{0}^{t}\int_{\mathbb{R}^{d}}{}\times\biggl[}{}
+H ( 2H-1 ) \int_{0}^{s} \bigl( Q (
B_{s-r}^{x},x )\nonumber\\
&&\hspace*{155pt}{}
-Q ( x,x ) \bigr ) r^{2H-2}\,dr\biggr]\,ds\,dx.
\nonumber
\end{eqnarray}
The convergence in $L^{2}$ to zero of $V_{\varepsilon}^{2}$ follows
from (%
\ref{Trace2}) and (\ref{Trace1}).
\end{pf}

\section{Skorohod type equation and Chaos expansion}\label{sec6}

In this section, we consider the following heat equation on $\mathbb{R}^{d}$
\begin{equation}
\cases{\displaystyle
\frac{\partial u}{\partial t}=\frac{1}{2}\Delta u+u\diamond\,\frac
{\partial W%
}{\partial t}(t,x) ,&\quad $t\geq0 , x\in\mathbb{R}^{d}$,\vspace*{2pt}\cr\displaystyle
u(0,x)=u_{0}(x) .%
}
\label{Wickequ}
\end{equation}
The difference between the above equation and equation (\ref{she}) is that
here we use the Wick product $\diamond$. This equation is studied in
Hu and
Nualart \cite{H-N07} for the case $H_{1}=\cdots=H_{d}=\frac{1}{2}$,
and in
\cite{hunuso} for the case $H_{1},\ldots ,H_{d}\in(\frac{1}{2},1) $,
$%
2H_{0}+H_{1}+\cdots+H_{d}>d+1.$ As in that paper, we can define the
following notion of solution.

\begin{definition}
An adapted random field $u=\{u ( t,x ) ,t\geq0 , x\in
\mathbb{R}%
^{d}\}$ such that $E ( u^{2} ( t,x )  ) <\infty$
for all $%
 ( t,x ) $ is a (mild) solution to equation (\ref
{Wickequ}) if for
any $ ( t,x ) \in\lbrack0,\infty)\times\mathbb{R}^{d}$, the
process $ \{ p_{t-s} ( x-y ) u ( s,y )
\mathbf{1}%
_{[0,t]} ( s )$, $s\geq0$, $y\in\mathbb{R}^{d} \} $ is
Skorohod integrable, and the following equation holds%
%
\begin{equation}
u ( t,x ) =p_{t}f ( x ) +\int_{0}^{t}\int
_{\mathbb{R}%
^{d}}p_{t-s} ( x-y ) u ( s,y ) \delta W_{s,y},
\label{def wick solu}
\end{equation}
where $p_{t} ( x ) $ denotes the heat kernel and
$p_{t}f (
x ) =\int_{\mathbb{R}^{d}}p_{t} ( x-y ) f (
y )\,dy$.
\end{definition}

From \cite{H-N07}, we know that the solution to equation (\ref{Wickequ})
exists with an explicit Wiener chaos expansion if and only if the Wiener
chaos expansion converges. Note that $g_{t,x}^{B} ( r,z )
:=\mathbf{%
1}_{[0,t]}(r)\mathbf{1}_{[0,B_{t-r}^{x}]}(z)\in\mathcal{H}$.
Formally, we
can write $g_{t,x}^{B} ( r,z ) =\delta (
B_{t-r}^{x}-z ) $
and we have%
\[
\int_{0}^{t}W ( dr,B_{s-r}^{x} ) =W (
g_{t,x}^{B} )
=\int_{0}^{t}\int_{\mathbb{R}^{d}}\delta ( B_{t-r}^{x}-z
) W (
dr,z )\,dz.
\]
Then in the same way as in Section 8 in \cite{hunuso} we can check
that $%
u ( t,x ) $ given by (\ref{Wick_solu}) below has the suitable
Wiener chaos expansion, which has to be convergent because $u (
t,x ) $ is square integrable. We state it as the following theorem.

\begin{theorem}
Suppose $H>\frac{1}{2}-\frac{1}{4}\gamma$ and $u_{0}$ is a bounded
measurable function. Then the unique (mild) solution to equation (\ref{Wickequ}) is given by the process
%
\begin{equation}
u ( t,x ) =E^{B} \bigl[ u_{0} ( B_{t}^{x} ) \exp
\bigl(W (
g_{t,x}^{B} ) -\tfrac{1}{2} \|
g_{t,x}^{B} \|_{\mathcal{H}}^{2}\bigr) \bigr] . \label{Wick_solu}
\end{equation}
\end{theorem}

\begin{remark}
We can also obtain a Feynman--Kac formula for the coefficients of the chaos
expansion of the solution to equation (\ref{she})%
\[
u ( t,x ) =\sum_{n=0}^{\infty}\frac{1}{n!}I_{n} (
h_{n} (
t,x )  )
\]
with
\[
h_{n} ( t,x ) =E^{B} \bigl[ u_{0} ( B_{t}^{x} )
g_{t,x}^{B} ( r_{1},z_{1} ) \cdots g_{t,x}^{B} (
r_{n},z_{n} ) \exp \bigl( \tfrac{1}{2} \|
g_{t,x}^{B} \|_{\mathcal{H}}^{2} \bigr)
\bigr] .
\]
\end{remark}

\begin{appendix}\label{appm}
\setcounter{equation}{0}
\setcounter{theorem}{0}
\section*{Appendix}

In this section, we denote by $B^{H}=\{B_{t}^{H},t\in\mathbb{R}\}$ a mean
zero Gaussian process with covariance $E(B_{t}^{H}B_{s}^{H})=
\frac{1}{2} (  \vert t \vert^{2H}+ \vert s
\vert
^{2H}-|t-s|^{2H} ) .$ Denote by $\mathcal{E}$ the space of all step
functions on $[-T,T]$. On $\mathcal{E}$, we introduce the following scalar
product $\langle\mathbf{1}_{[0,t]},\mathbf{1}_{[0,s]}\rangle
_{\mathcal{H}%
_{0}}=R_{H}(t,s),$ where if $t<0$ we assume that $\mathbf{1}_{[0,t]}=-%
\mathbf{1}_{[t,0]}$. Let $\mathcal{H}_{0}$ be the closure of
$\mathcal{E}$
with respect to the above scalar product.\vadjust{\goodbreak}

For $r>0$, $\varepsilon>0$ and $\beta>0,$ let%
\[
f^{\varepsilon} ( r ) :=\frac{1}{4\varepsilon^{2}} [
2r^{\beta}- \vert r-2\varepsilon \vert^{\beta}- (
r+2\varepsilon ) ^{\beta} ] .
\]
It is easy to see that
%
\begin{equation}\label{DCT_lmt}
\lim_{\varepsilon\downarrow0}f^{\varepsilon} ( r )
=\beta (
\beta-1 ) r^{\beta-2}.
\end{equation}

%
%
%

\begin{lemma}
For any $r>0$, $\varepsilon>0$ and $0<\beta<2,$
%
\begin{equation}
 \vert f^{\varepsilon} ( r )  \vert\leq
64r^{\beta-2}.
\label{DCT_ctrl}
\end{equation}
\end{lemma}

\begin{pf}
If $0<r<4\varepsilon,$ then $ \vert r-2\varepsilon \vert
^{\beta
}< ( 2\varepsilon ) ^{\beta}$, $ ( r+2\varepsilon
 )
^{\beta}< ( 6\varepsilon ) ^{\beta}$, and hence (noting
that $\beta<2$)
\[
 \vert f^{\varepsilon} ( r )  \vert\leq
4^{\beta
+1}\varepsilon^{\beta-2}\leq64r^{\beta-2}.
\]
On the other hand, if $r\geq4\varepsilon,$ then
\begin{eqnarray*}
r^{\beta}- \vert r-2\varepsilon \vert^{\beta}
&=&2\varepsilon
\beta\int_{0}^{1} ( r-2\lambda\varepsilon ) ^{\beta
-1}\,d\lambda,
\\
r^{\beta}- ( r+2\varepsilon ) ^{\beta} &=&-2\varepsilon
\beta
\int_{0}^{1} ( r+2\lambda\varepsilon ) ^{\beta
-1}\,d\lambda,
\end{eqnarray*}
and hence
\begin{eqnarray*}
f^{\varepsilon} ( r ) &=&\frac{1}{2\varepsilon}\beta\int
_{0}^{1}%
 [  ( r-2\lambda\varepsilon ) ^{\beta-1}- (
r+2\lambda
\varepsilon ) ^{\beta-1} ]\,d\lambda\\
&=&2\beta ( \beta-1 ) \int_{0}^{1}\int_{0}^{1}\lambda
 (
r-2\lambda\varepsilon+4\mu\lambda\varepsilon ) ^{\beta
-2}\,d\mu\,d\lambda.
\end{eqnarray*}
Therefore, using $\beta<2$ and $r\geq4\varepsilon$ we obtain
\[
 \vert f^{\varepsilon} ( r )  \vert\leq2\beta
 (
r-2\varepsilon ) ^{\beta-2}\leq4r^{\beta-2} \biggl( \frac{%
r-2\varepsilon}{r} \biggr) ^{\beta-2}\leq16r^{\beta-2}.
\]
\upqed
\end{pf}

\begin{lemma}
\label{DCT Lmt phi}For any $s>0$, $0<\beta<1$ and $\phi\in C^{\alpha
} (  [ 0,T ]  ) $ with $\alpha>1-\beta$,
%
\begin{equation}
  \lim_{\varepsilon\rightarrow0}\int_{0}^{s}\phi ( r )
f^{\varepsilon} ( r )\,dr=\phi ( 0 ) \beta
s^{\beta
-1}+\beta ( \beta-1 ) \int_{0}^{s} \bigl( \phi (
r )
-\phi ( 0 )  \bigr) r^{\beta-2}\,dr. \label{phi_Lmt}\hspace*{-35pt}
\end{equation}
Moreover,
%
\begin{equation}
 \biggl\vert\int_{0}^{s}\phi ( r ) f^{\varepsilon} (
r )\,dr \biggr\vert\leq C ( \beta,\alpha )  (  \Vert
\phi
 \Vert_{\infty}s^{\beta-1}+ \Vert\phi \Vert
_{\alpha
}s^{\alpha+\beta-1} ) . \label{phi_ctrl}
\end{equation}
\end{lemma}

\begin{pf}
The lemma follows easily from (\ref{iprdct1}) and (\ref{DCT_ctrl}) if we
rewrite
\[
\int_{0}^{s}\phi ( r ) f^{\varepsilon} ( r )\,dr=\phi
 ( 0 ) \int_{0}^{s}f^{\varepsilon} ( r )\,dr+\int_{0}^{s}
 [ \phi ( r ) -\phi ( 0 )  ]
f^{\varepsilon
} ( r )\,dr.
\]%
\upqed
\end{pf}\eject

\begin{lemma}
\label{inner_prdct1}For any bounded function $\phi\in\mathcal
{H}_{0}$ and
any $s,t\geq0$, we have
%
\begin{equation}
 \bigl \langle\mathbf{1}_{[0,s]}\phi,\mathbf{1}_{[0,t]}
\bigr\rangle_{%
\mathcal{H}_{0}}=H\int_{0}^{s}\phi ( r )  [
r^{2H-1}+\operatorname{sign} ( t-r )  \vert t-r \vert^{2H-1} ]\,dr.
\label{iprdct1}\hspace*{-35pt}
\end{equation}
If $u<s<t$, we have
%
\begin{equation}
  \bigl\langle\mathbf{1}_{[0,s]}\phi,\mathbf{1}_{[u,t]}
\bigr\rangle_{%
\mathcal{H}_{0}}=H\int_{0}^{s}\phi ( r )  [  (
t-r )
^{2H-1}-\operatorname{sign} ( u-r )  \vert u-r \vert
^{2H-1}%
 ]\,dr.\label{iprdct2}\hspace*{-35pt}
\end{equation}
\end{lemma}

\begin{pf}
We only have to prove (\ref{iprdct1}) since (\ref{iprdct2}) follows easily.
Without loss of generality, assume that $\phi=\sum
_{i=1}^{n}a_{i}\mathbf{%
1}_{ [ t_{i-1},t_{i} ] }$, where $0=t_{0}\leq t_{1}\leq
\cdots\leq
t_{n}=s$. (If $t<s$, we assume that $t=t_{i}$ for some $0<i<n.$) Then%
\begin{eqnarray*}
 \bigl \langle\mathbf{1}_{[0,s]}\phi,\mathbf{1}_{[0,t]}
\bigr\rangle_{%
\mathcal{H}_{0}}&=&E\sum_{i=1}^{n}a_{i} (
B_{t_{i}}^{H}-B_{t_{i-1}}^{H} ) B_{t}^{H} \\
&=&\sum_{i=1}^{n}a_{i}\frac{1}{2} (
t_{i}^{2H}-t_{i-1}^{2H}+ \vert t-t_{i-1} \vert^{2H}-
\vert
t-t_{i} \vert^{2H} ) \\
&=&H\int_{0}^{s}\phi ( r )  [ r^{2H-1}+\operatorname
{sign} (
t-r )  \vert t-r \vert^{2H-1} ]\,dr.
\end{eqnarray*}
\upqed
\end{pf}

Using Lemma \ref{inner_prdct1} and similar arguments to those in the proof
of Lem\-ma~\ref{DCT Lmt phi}, we can prove the following lemma.
%
%

\begin{lemma}
\label{DCT_Lm2}For any $s>0$, for any $\phi\in C^{\alpha} ( %
 [ 0,T ]  ) $ with $\alpha>1-2H$,
\begin{eqnarray*}
\lim_{\varepsilon\rightarrow0} \biggl\langle\mathbf{1}_{[0,s]}\phi
,\frac{%
1}{2\varepsilon}\mathbf{1}_{[s-\varepsilon,s+\varepsilon]}
\biggr\rangle_{%
\mathcal{H}_{0}}
=\phi ( s ) Hs^{2H-1}+c_0 \int_{0}^{s}\bigl (
\phi ( s-r ) -\phi ( s )  \bigr) r^{2H-2}\,dr,
\end{eqnarray*}
where $c_0=H(2H-1)$. Moreover,
%
\begin{equation}
 \biggl| \biggl\langle\mathbf{1}_{[0,s]}\phi,\frac
{1}{2\varepsilon}%
\mathbf{1}_{[s-\varepsilon,s+\varepsilon]}\biggr \rangle_{\mathcal
{H}%
_{0}} \biggr|\leq C ( H,\alpha )  (  \|
\phi \|_{\infty}s^{2H-1}+ \|\phi
 \|_{\alpha}s^{\alpha+2H-1} ) .
\label{DCT2_ctrl}\hspace*{-35pt}
\end{equation}
\end{lemma}

\begin{pf}
Applying Lemma \ref{inner_prdct1} and making a substitution, we get
\begin{eqnarray*}
&& \biggl\langle\mathbf{1}_{[0,s]}\phi,\frac{1}{2\varepsilon
}\mathbf{1}%
_{[s-\varepsilon,s+\varepsilon]} \biggl\rangle_{\mathcal{H}_{0}} \\
&& \qquad =\frac{H}{2\varepsilon}\int_{0}^{s}\phi ( s-u )
[  (
u+\varepsilon ) ^{2H-1}-\operatorname{sign} ( u-\varepsilon
 )
 \vert u-\varepsilon \vert^{2H-1} ]\,du \\
&& \qquad =:H\phi ( s ) \int_{0}^{s}g^{\varepsilon} (
u )\,du+H\int_{0}^{s} [ \phi ( s-u ) -\phi ( s
)  ]
g^{\varepsilon} ( u )\,du,
\end{eqnarray*}
where we let
\[
g^{\varepsilon} ( u ) =\frac{1}{2\varepsilon} [
 (
u+\varepsilon ) ^{2H-1}-\operatorname{sign} ( u-\varepsilon
 )
 \vert u-\varepsilon \vert^{2H-1} ] .\vadjust{\goodbreak}
\]
If $0<u <2\varepsilon$, we have $ \vert g^{\varepsilon} (
u )  \vert\leq16r^{2H-2}.$ On the other hand, if
$u>2\varepsilon
$,
\begin{eqnarray*}
 \vert g^{\varepsilon} ( u )  \vert&=&
\biggl\vert\frac{1%
}{2\varepsilon} [  ( u-\varepsilon ) ^{2H-1}- (
u+\varepsilon ) ^{2H-1} ]  \biggr\vert\\
&=&\frac{1}{2} ( 1-2H ) \int_{-1}^{1} ( u-\lambda
\varepsilon
 ) ^{2H-2}\,d\lambda\leq ( 1-2H ) u^{2H-2}.
\end{eqnarray*}
Then the lemma follows by noticing that $\lim_{\varepsilon\rightarrow
0}g^{\varepsilon} ( u ) = ( 2H-1 ) u^{2H-2}.$
\end{pf}

%

\begin{lemma}
\label{DCT_Lm4}
For any $\phi\in C^{\alpha} (  [ 0,T ]
 ) $ with $\alpha>1-2H$, for any $s>0$,
%
\begin{eqnarray}
\label{DCT4_Lmt}
&&\lim_{\varepsilon\rightarrow0} \biggl\langle\frac{1}{2\varepsilon}
\int_{0}^{s}\mathbf{1}_{[\theta-\varepsilon,\theta+\varepsilon
]}\phi
 ( \theta )\,d\theta,\frac{1}{2\varepsilon}\mathbf{1}%
_{[s-\varepsilon,s+\varepsilon]} \biggr\rangle_{\mathcal{H}_{0}} \nonumber\hspace*{-35pt}
\\[-8pt]
\\[-8pt]
&& \qquad =\phi ( s ) Hs^{2H-1}+H ( 2H-1 ) \int
_{0}^{s} \bigl(
\phi ( s-r ) -\phi ( s )  \bigr) r^{2H-2}\,dr.
\nonumber\hspace*{-35pt}
\end{eqnarray}
Moreover,
%
\begin{eqnarray}
\label{DCT4_ctrl}
&&\biggl \vert \biggl\langle\frac{1}{2\varepsilon}\int_{0}^{s}\mathbf
{1}%
_{[\theta-\varepsilon,\theta+\varepsilon]}\phi ( \theta
 )\,d\theta,\frac{1}{2\varepsilon}\mathbf{1}_{[s-\varepsilon
,s+\varepsilon
]} \biggl\rangle_{\mathcal{H}_{0}} \biggl\vert \nonumber
\\[-8pt]
\\[-8pt]
&& \qquad \leq C ( H,\alpha )  (  \|\phi
 \|_{\infty}Hs^{2H-1}+ \|
\phi
 \|_{\alpha}s^{\alpha+2H-1} ) .
\nonumber
\end{eqnarray}
\end{lemma}

\begin{pf}
By Fubini's theorem and making a substitution, we have%
\begin{eqnarray*}
&& \biggl\langle\frac{1}{2\varepsilon}\int_{0}^{s}\mathbf
{1}_{[\theta
-\varepsilon,\theta+\varepsilon]}\phi ( \theta )\,d\theta,
\frac{1}{2\varepsilon}\mathbf{1}_{[s-\varepsilon,s+\varepsilon
]} \biggl\rangle_{\mathcal{H}_{0}} \\
&& \qquad =\frac{1}{4\varepsilon^{2}}E \biggl[ \int_{0}^{s}\phi (
\theta )
 ( B_{\theta+\varepsilon}^{H}-B_{\theta-\varepsilon}^{H} )
 ( B_{s+\varepsilon}^{H}-B_{s-\varepsilon}^{H} )\,d\theta
 \biggr]
\\
&& \qquad =\frac{1}{8\varepsilon^{2}}\int_{0}^{s}\phi ( s-\theta
 )  [
2r^{2H}- \vert r-2\varepsilon \vert^{2H}- (
r+2\varepsilon
 ) ^{2H} ]\,dr.
\end{eqnarray*}
Then (\ref{DCT4_Lmt}) and (\ref{DCT4_ctrl}) follow from Lemma \ref
{DCT Lmt
phi}.
\end{pf}
\end{appendix}


%

\printaddresses


\begin{thebibliography}{14}

\bibitem{BC95}
%
\begin{barticle}[mr]
\bauthor{\bsnm{Bertini},~\bfnm{Lorenzo}\binits{L.}} \AND
\bauthor{\bsnm{Cancrini},~\bfnm{Nicoletta}\binits{N.}}
(\byear{1995}).
\btitle{The stochastic heat equation: {F}eynman--{K}ac formula and
intermittence}.
\bjournal{J. Statist. Phys.}
\bvolume{78}
\bpages{1377--1401}.
\bid{issn={0022-4715}, mr={1316109}}
\end{barticle}
%
\endbibitem

\bibitem{HJT}
%
\begin{bmisc}[auto:STB|2011-03-03|12:04:44]
\bauthor{\bsnm{Hu},~\bfnm{Y.}\binits{Y.}},
\bauthor{\bsnm{Jolis},~\bfnm{M.}\binits{M.}} \AND
\bauthor{\bsnm{Tindel},~\bfnm{S.}\binits{S.}}
(\byear{2011}).
\bhowpublished{On Stratonovich and Skorohod stochastic calculus for Gaussian
processes.  Preprint. Available at
\texttt{\href{http://arxiv.org/abs/1101.3441}{http://arxiv.org/abs/}
\href{http://arxiv.org/abs/1101.3441}{1101.3441}}.}
\end{bmisc}
%
\endbibitem

\bibitem{H-N07}
%
\begin{barticle}[mr]
\bauthor{\bsnm{Hu},~\bfnm{Yaozhong}\binits{Y.}} \AND
\bauthor{\bsnm{Nualart},~\bfnm{David}\binits{D.}}
(\byear{2009}).
\btitle{Stochastic heat equation driven by fractional noise and local time}.
\bjournal{Probab. Theory Related Fields}
\bvolume{143}
\bpages{285--328}.
\bid{doi={10.1007/s00440-007-0127-5}, issn={0178-8051}, mr={2449130}}
\end{barticle}
%
\endbibitem

\bibitem{hunuso}
%
\begin{barticle}[auto:STB|2011-03-03|12:04:44]
\bauthor{\bsnm{Hu},~\bfnm{Y.}\binits{Y.}},
\bauthor{\bsnm{Nualart},~\bfnm{D.}\binits{D.}} \AND
\bauthor{\bsnm{Song},~\bfnm{J.}\binits{J.}}
(\byear{2011}).
\btitle{Feynman--Kac formula for heat equation driven by fractional white
noise}.
\bjournal{Ann. Probab.}
\bvolume{39}
\bpages{291--326}.
\end{barticle}
%
\endbibitem

\bibitem{KR}
%
\begin{bmisc}[auto:STB|2011-03-03|12:04:44]
\bauthor{\bsnm{Kruk},~\bfnm{I.}\binits{I.}} \AND
\bauthor{\bsnm{Russo},~\bfnm{F.}\binits{F.}}
(\byear{2010}).
\bhowpublished{Malliavin--Skorohod calculus and Paley--Wiener integral for
covariance singular processes.  Preprint. Available at
\texttt{\href{http://arxiv.org/abs/1011.6478}{http://arxiv.org/}
\href{http://arxiv.org/abs/1011.6478}{abs/1011.6478}}}.
\end{bmisc}
%
\endbibitem

\bibitem{M-V05}
%
\begin{barticle}[mr]
\bauthor{\bsnm{Mocioalca},~\bfnm{Oana}\binits{O.}} \AND
\bauthor{\bsnm{Viens},~\bfnm{Frederi}\binits{F.}}
(\byear{2005}).
\btitle{Skorohod integration and stochastic calculus beyond the fractional
{B}rownian scale}.
\bjournal{J. Funct. Anal.}
\bvolume{222}
\bpages{385--434}.
\bid{doi={10.1016/j.jfa.2004.07.013}, issn={0022-1236}, mr={2132395}}
\end{barticle}
%
\endbibitem

\bibitem{Nu06}
%
\begin{bbook}[mr]
\bauthor{\bsnm{Nualart},~\bfnm{David}\binits{D.}}
(\byear{2006}).
\btitle{The {M}alliavin Calculus and Related Topics},
\bedition{2nd} ed.
\bpublisher{Springer}, \baddress{Berlin}.
\bid{mr={2200233}}
\end{bbook}
%
\endbibitem

\bibitem{Nu-Sch01}
%
\begin{barticle}[mr]
\bauthor{\bsnm{Nualart},~\bfnm{David}\binits{D.}} \AND
\bauthor{\bsnm{Schoutens},~\bfnm{Wim}\binits{W.}}
(\byear{2001}).
\btitle{Backward stochastic differential equations and {F}eynman--{K}ac formula
for {L}\'evy processes, with applications in finance}.
\bjournal{Bernoulli}
\bvolume{7}
\bpages{761--776}.
\bid{doi={10.2307/3318541}, issn={1350-7265}, mr={1867081}}
\end{barticle}
%
\endbibitem

\bibitem{O-P93}
%
\begin{barticle}[mr]
\bauthor{\bsnm{Ocone},~\bfnm{Daniel}\binits{D.}} \AND
\bauthor{\bsnm{Pardoux},~\bfnm{{\'E}tienne}\binits{{\'E}.}}
(\byear{1993}).
\btitle{A stochastic {F}eynman--{K}ac formula for anticipating {SPDE}s, and
application to nonlinear smoothing}.
\bjournal{Stochastics Stochastics Rep.}
\bvolume{45}
\bpages{79--126}.
\bid{issn={1045-1129}, mr={1277363}}
\end{barticle}
%
\endbibitem

\bibitem{OS02}
%
\begin{barticle}[mr]
\bauthor{\bsnm{Ouerdiane},~\bfnm{Habib}\binits{H.}} \AND
\bauthor{\bsnm{Silva},~\bfnm{Jos{\'e}~Luis}\binits{J.~L.}}
(\byear{2002}).
\btitle{Generalized {F}eynman--{K}ac formula with stochastic potential}.
\bjournal{Infin. Dimens. Anal. Quantum Probab. Relat. Top.}
\bvolume{5}
\bpages{243--255}.
\bid{doi={10.1142/S0219025702000808}, issn={0219-0257}, mr={1914836}}
\end{barticle}
%
\endbibitem

\bibitem{SKM}
%
\begin{bbook}[mr]
\bauthor{\bsnm{Samko},~\bfnm{Stefan~G.}\binits{S.~G.}},
\bauthor{\bsnm{Kilbas},~\bfnm{Anatoly~A.}\binits{A.~A.}} \AND
\bauthor{\bsnm{Marichev},~\bfnm{Oleg~I.}\binits{O.~I.}}
(\byear{1993}).
\btitle{Fractional Integrals and Derivatives}.
\bpublisher{Gordon and Breach Science Publishers}, \baddress{Yverdon}.
\bid{mr={1347689}}
\end{bbook}
%
\endbibitem


\bibitem{Zahle98}
%
\begin{barticle}[mr]
\bauthor{\bsnm{Z{\"a}hle},~\bfnm{M.}\binits{M.}}
(\byear{1998}).
\btitle{Integration with respect to fractal functions and stochastic calculus.
{I}}.
\bjournal{Probab. Theory Related Fields}
\bvolume{111}
\bpages{333--374}.
\bid{doi={10.1007/s004400050171}, issn={0178-8051}, mr={1640795}}
\end{barticle}
%
\endbibitem

\end{thebibliography}
\end{document}